\documentclass{amsart}

\usepackage{amssymb}
\usepackage{mathrsfs}

\usepackage[cmtip,all]{xy}

\usepackage{graphicx}
\usepackage{bbm}
\usepackage{xcolor}
\usepackage{stmaryrd}
\copyrightinfo{2022}{American Mathematical Society}

\newtheorem{theorem}{Theorem}[section]

\theoremstyle{definition}

\theoremstyle{remark}
\newtheorem{remark}[theorem]{Remark}

\numberwithin{equation}{section}

\newcommand{\CC}{\mathbb{C}}
\newcommand{\QQ}{\mathbb{Q}}

\newcommand{\FF}{\mathbb{F}}

\newcommand{\Sym}{\operatorname{Sym}}

\newcommand{\rs}{{\operatorname{rs}}}

\newcommand{\Gm}{\mathbb{G}_m}
\newcommand{\Ga}{\mathbb{G}_a}
\newcommand{\GL}{\operatorname{GL}}
\newcommand{\Mat}{\operatorname{Mat}}
\newcommand{\PGL}{\operatorname{PGL}}
\newcommand{\SL}{\operatorname{SL}}
\newcommand{\Sp}{\operatorname{Sp}}

\newcommand{\SO}{{\operatorname{SO}}}
\newcommand{\GO}{{\operatorname{GO}}}

\newcommand{\Gal}{\operatorname{Gal}}

\newcommand{\tr}{\operatorname{tr}}
\newcommand{\Spec}{\operatorname{Spec\,}}
\newcommand{\Vol}{\operatorname{Vol}}
\newcommand{\diag}{{\operatorname{diag}}}

\newcommand{\reg}{{\operatorname{reg}}}

\newcommand{\rk}{{\operatorname{rk}}}

\newcommand{\Std}{{\operatorname{Std}}}
\newcommand{\st}{{\operatorname{st}}}
\newcommand{\Ad}{{\operatorname{Ad}}}

\newcommand{\Res}{\operatorname{Res}}

\newcommand{\TF}{{\operatorname{TF}}}
\newcommand{\STF}{{\operatorname{STF}}}
\newcommand{\RTF}{{\operatorname{RTF}}}
\newcommand{\KTF}{{\operatorname{KTF}}}

\newcommand{\LG}{{^LG}}
\newcommand{\LH}{{^LH}}

\newcommand{\Dfrac}[2]{%
  \ooalign{%
    $\genfrac{}{}{1.2pt}0{#1}{#2}$\cr%
    $\color{white}\genfrac{}{}{.4pt}0{\phantom{#1}}{\phantom{#2}}$}%
}

\begin{document}

\title[Local and global questions ``beyond endoscopy'']{Local and global questions ``beyond endoscopy''}

\author{Yiannis Sakellaridis}
\address{Department of Mathematics, Johns Hopkins University, Baltimore, MD 21218, USA.}
\email{sakellar@jhu.edu}
\thanks{}


\begin{abstract}
 The near-completion of the program of endoscopy poses the question of what lies next. 
 This article takes a broad view of ideas beyond the program of endoscopy, highlighting the connections among them, and emphasizing the relationship between local and global aspects. Central among those ideas is the one proposed in a 2000 lecture of R.~P.~Langlands, aiming to extract from the stable trace formula of a group $G$ the bulk of those automorphic representations in the image of the conjectural functorial lift corresponding to a morphism of $L$-groups $\LH\to\LG$. With the extension of the problem of functoriality to the ``relative'' setting of spherical varieties and related spaces, some structure behind such comparisons has started to reveal itself. In a seemingly unrelated direction, a program initiated by Braverman--Kazhdan, also around 2000, to generalize the Godement--Jacquet proof of the functional equation to arbitrary $L$-functions, has received renewed attention in recent years. We survey ideas and developments in this direction, as well, and discuss the relationship between the two programs.
\end{abstract}

\maketitle

\tableofcontents

\section{Introduction}

\subsection{The Langlands program and $L$-functions}

The desire to understand $L$-functions animated the birth of the Langlands program, as is evident from the 1967 letter of Robert P.~Langlands to Andr\'e Weil \cite{Langlands-Weil}. Once in place, however, the conjectures of reciprocity and functoriality seem to throw the direct study of $L$-functions out of the window; in the creator's own words, ``since functoriality, once established, entails immediately the analytic continuation of all the functions $L(s, \pi , \rho)$, [..] there would seem to be little point in pursuing seriously the various methods introduced for dealing with this or that special case'' \cite{Langlands-wherestands}, that is, the methods that one uses to represent automorphic $L$-functions. 

Fortunately, such methods continued to be the object of intense study by mathematicians such as Piatetski-Shapiro, Jacquet, and Shalika, despite the fact that ``the mathematical community was not
very supportive because the project deviated from
the accepted canon and was regarded as perverse'' \cite{JaNotices}.

Then, in a strange twist, $L$-functions came again to play a prominent role in Langlands' vision of the path to functoriality ``Beyond Endoscopy'' \cite{Langlands-BE}: It was no more sufficient to work with the trace formula as developed by Selberg and Arthur, out of ``pure'' (albeit extremely complicated) ingredients such as compactly supported test functions and their orbital integrals. One also had to introduce $L$-functions ``by hand,'' whose poles had the role of picking up the images of functorial lifts. 

``By hand'' could mean writing down the Dirichlet series of an automorphic $L$-function 
\[ L(s, \pi , \rho) = \sum_{n\ge 1} \frac{a_n}{n^s},\]
where the $a_n$'s, for prime powers $n=p^i$, are polynomials in the Satake parameters (or: Frobenius--Hecke parameters, as Langlands prefers to call them) of the automorphic representation $\pi$ at $p$. Inverting the Satake isomorphism (and glossing over ramified places, for the purposes of this introduction), we obtain a corresponding series $h_{r,s}$ of Hecke operators, whose trace on the automorphic representation $\pi$ is the $L$-value above. The insertion of $h_{r,s}$ as a ``nonstandard test function'' into the trace formula gives, at least formally, a spectral expansion where the usual term for $\pi$ is weighed by the factor $L(s, \pi , \rho)$.

The idea of Langlands (presented and expanded by Arthur in \cite{Arthur-BE}) was that, weighing the trace formula by $L$-functions, and taking residues at the poles of those, one could detect the images of functorial lifts. The original idea, actually, involved weighing by the \emph{logarithmic derivatives} of $L$-functions, so that the residues would be weighed only by the orders of the poles; it is generally understood, however, that the geometric expressions for the logarithmic derivative are prohibitively difficult to manage.

The Beyond Endoscopy project has seen some limited successes, in the two decades since its inception; to mention a few: In his thesis \cite{Venkatesh}, which used the Kuznetsov, rather than the Selberg trace formula, Venkatesh  successfully employed ``Beyond Endoscopy'' techniques to re-establish functoriality from tori to $\GL_2$. A series of papers of Herman \cite{Herman-RS, Herman} used similar techniques to give a new proof of analytic continuation and the functional equation of the standard and Rankin--Selberg $L$-functions for $\GL_2$. 
Altu\u{g}'s work \cite{Altug1,Altug2,Altug3} managed to rigorously complete the argument outlined by Langlands for $\GL_2$, isolating the contribution of ``special representations'' from the Arthur--Selberg trace formula, and showing that the standard $L$-function does not ``pick up'' functorial lifts in the cuspidal spectrum. The scarcity of such successes attests the impressive nature of these feats, and the abilities of the authors. At first sight, however, these works employ case-specific brute-force techniques in order to manage the difficult analytic expressions that arise after inserting $L$-functions into the trace formula, and do not necessarily point to a paradigm that can be scaled up.

The present article is predicated on the conviction that the insertion of $L$-functions into the trace formula gives rise to structures that need to be understood using tools from representation theory and geometry. This is grounded on the limited experience that I have obtained through my own involvement in the project, and the difficulty of reaching concrete conclusions based simply on techniques of analytic number theory. In particular, we will attempt to mimic the paradigm of endoscopy, where local results and comparisons are formulated, and fed into the trace formula to produce the desired global results. 

For example, as was understood by Ng\^o \cite{Ngo-PS}, following his joint work with Frenkel \cite{FN}, a ``geometrization'' of the Euler factors of the aforementioned series $h_{r,s}$ of Hecke elements is possible by replacing the reductive group $G$ by a suitable reductive monoid $G'\supset G$ (possibly singular), generalizing the Godement--Jacquet theory for $\Mat_n\supset \GL_n$. The Hecke series (for a local field of equal characteristic) is then the Frobenius trace on an intersection complex associated to the arc space of this monoid \cite{BNS}. This makes the above trace formula, weighed by $L$-functions, amenable, at least in principle, to the geometric techniques and insights that came with the proof of the Fundamental Lemma.

This generalization of Godement--Jacquet theory, including the suggestion to look at reductive monoids, had been the topic of a series of papers by Braverman and Kazhdan, also around the year 2000 \cite{BK2,BK4}. Moreover, reductive monoids are but a special case of affine spherical varieties, whose putative Schwartz spaces, as noted in \cite{SaRS},  should explain and generalize a large portion of the Rankin--Selberg method of integral representations of $L$-functions. 

Broadening the problem of functoriality to more general spaces than reductive groups, and admitting the relative trace formula of Jacquet into our purview of desired ``Beyond Endoscopy'' comparisons, gives us a broader playing ground for testing the core ideas. In particular, a \emph{local} theory of ``Beyond Endoscopy'' comparisons has been developed by the author in low-rank cases \cite{SaBE1, SaBE2, SaTransfer1, SaTransfer2, SaRankone}, revealing surprising uniformity of the pertinent ``transfer operators,'' and suggesting structures that might scale up to higher rank. These structures seem to be controlled by $L$-functions. Thus, paradoxically, in order to address functoriality ``beyond endoscopy,'' one is led to deviate from
the canon in a rather systematic way.

\subsection{Outline of this article}

``Beyond Endoscopy'' (capitalized when referring to the original proposal, lowercase for a broader set of ideas) was never a proposal with the specificity and solid evidence of endoscopy, that Langlands himself provided, along with Labesse, Shelstad, and Kottwitz, in the late 70s and 80s \cite{LL, LS, Kottwitz-rational}. Its development heretofore took a toll on the careers of a few young and promising mathematicians, and we still do not have a broadly accepted view of its short- and medium-term goals, and its accomplishments to date. 

Thus, it would be futile to attempt an all-encompassing overview of the subject, that would also be mathematically meaningful. The present article takes an ahistorical approach, trying to piece together several concrete problems, strategies, and results developed beyond the purview of endoscopy, and to highlight their connections, whether they were developed with Langlands' proposal in mind, or not. In doing so, it does not do justice to the entire body of work that was and continues to be developed under the ``Beyond Endoscopy'' label, and which has produced useful insights into various problems of analytic or arithmetic nature. This account is certainly colored by my own limitations and my personal vision, and will not present everyone's point of view of Beyond Endoscopy.

We start by discussing local aspects. After introducing transfer operators as the main vehicle of Beyond Endoscopy, our point of departure becomes Langlands' ``Singularit\'es et transfert'' \cite{Langlands-ST}, where those operators are studied for the functorial lift between tori and the stable trace formula of $\SL_2$ (or $\GL_2$). Section \ref{sec:localtori} concludes with a different way to realize this local transfer, replacing the trace formula by the Kuznetsov formula.

In Section \ref{sec:samerank} we continue with the local theory, discussing the problem of stable transfer for identical $L$-groups. This problem only makes sense in the broader setting of the relative trace formula, where different spherical varieties can share the same $L$-group. A particular case is the transfer between the stable and Kuznetsov formulas for $\SL_2$ (or $\GL_2$), which, I argue later, lies behind the global theory of Beyond Endoscopy developed by Altu\u{g}, as well as the thesis of Rudnick \cite{Rudnick}, that predated Beyond Endoscopy.

In Section \ref{sec:Schwartz} we take a look at various ``Schwartz spaces'' of functions that encode $L$-functions, while in Section \ref{sec:Hankel} we discuss ways to give meaning to the functional equations of $L$-functions, either by means of ``nonlinear Fourier transforms'' on appropriate spherical varieties, or by their analogs, that we call Hankel transforms, at the level of trace formulas. 

Finally, in Section \ref{sec:global} we discuss global comparisons of trace formulas ``Beyond Endoscopy,'' in the light of the preceding local discussion.

A particularly bewildering theme beyond endoscopy is the role of Fourier transforms and the Poisson summation formula, on spaces that have no obvious linear structure, such as the invariant-theoretic quotient of $G$ by conjugation (the space of characteristic polynomials, when $G=\GL_n$). Such a formula is completely paradoxical, since the space of conjugacy classes does not carry any natural linear structure. In \S~\ref{ss:transferrankone}, we will discuss symplectic structures which make the appearance of Fourier transforms more natural. The exact relationship between ``beyond endoscopy'' methods and symplectic geometry, however, is a largely unexplored topic, and poses a promising challenge for creating better foundations for the program in the near future.

\subsection{Notation}

\begin{itemize}
 \item For an affine $G$-variety $X$ over a field $F$, the invariant-theoretic quotient $\Spec F[X]^G$ will be denoted by $X\sslash G$. For $G$ acting on itself by conjugation, this quotient will be denoted $\Dfrac{G}{G}$.
 \item When the meaning is clear from the context, for a variety $X$ over a local field $F$, we denote the set $X(F)$ simply by $X$.
 \item For a variety $X$ over a local field $F$, we \emph{usually} denote by $\mathcal S(X)$ the space of Schwartz \emph{measures} on the $F$-points of $X$. This space makes sense when $X$ is smooth, in which case ``Schwartz measures'' means ``smooth of rapid decay'' in the Archimedean case (see \cite{AGSchwartz}), and ``smooth of compact support'' in the non-Archimedean case. At various points in the article, we speculate on extensions of this notion for non-smooth varieties. \emph{Sometimes} we use this notation for functions, instead of measures -- and we say so. Finally, the notation $\mathcal D(X)$ is used for Schwartz half-densities on $X(F)$. Often, $X(F)$ is homogeneous and admits an invariant positive measure $dx$ under a group $G(F)$, in which case one can easily pass from functions to half-densities and from half-densities to measures by multiplying by $(dx)^\frac{1}{2}$.
 \item For a group $G$ acting (on the left or right) on a space $X$, we denote by $\star$ the convolution action of measures on $G$ on functions/half-densities/measures on $X$. To be clear, convolution is the pushforward under the action map, i.e., 
 \[ (\mu\star f)(x) = \int_G f(x\cdot g^{-1}) \mu(g).\]
\end{itemize}



\subsection{Convention on citations}

Since a significant portion of work on the subject has appeared in lecture notes and non-refereed preprints, I will take a cautious approach in stating and attributing results. Although unpublished results and claims will be discussed, only published results, and their authors, appear with the label ``Theorem.''

\subsection{Acknowledgments}

I thank B.\ C.\ Ng\^o for useful feedback, and for sharing with me his notes on L.\ Lafforgue's Fourier kernel.

\section{Local stable transfer; the lift from tori to $\SL_2$.} \label{sec:localtori}

\subsection{Stable test measures, functions, and half-densities}

We start with the local theory of ``Beyond Endoscopy.'' Fix a local (locally compact) field $F$, and a quasisplit connected reductive group $G$ over $F$. Let $\mathcal S(G)$ be the space of Schwartz measures on $G=G(F)$, and denote the image of the pushforward map 
\[ \mathcal S(G) \xrightarrow{p} \operatorname{Meas}(\Dfrac{G}{G}),\]
by $\mathcal S(\frac{G}{G})^\st$, the space of \emph{stable} test measures (for the adjoint quotient of the group). Explicitly, by the Weyl integration formula, if $f=\Phi \, dg$ is a Schwartz measure, written as the product of a Schwartz function by a Haar measure, its pushforward $f'=p_! f$ to $\Dfrac{G}{G}$ reads, on the strongly regular semisimple set, 
\begin{equation}\label{pushfmeasure} f' ([t]) = \SO_t (\Phi) |D(t)| d[t],
\end{equation}
where
\begin{itemize}
 \item $[t]$ denotes the stable conjugacy class of a rational, strongly regular element $t \in T^\reg(F)$, for some maximal torus $T \subset G$;
 \item $d[t]$ denotes a measure on the image of $T^\reg(F)$ in $\Dfrac{G}{G} = T\sslash W$ whose pullback\footnote{By ``pullback'' of a Borel measure under a map $p:X\to Y$ which is a local homeomorphism, we will mean the measure on $X$ which, on every open $U\subset X$ which maps bijectively onto its image, coincides with the given measure on $p(U)$.} to $T^\reg$ coincides with some Haar measure $dt$ on $T(F)$;
 \item $\SO_t$ denotes the stable orbital integral of $\Phi$ over the stable orbit of $t$, defined with respect to the quotient measure $dg/dt$ on the conjugacy classes (which are isomorphic to $T(F)\backslash G(F)$);
 \item $D(t) = \det(I - \Ad_{\mathfrak g/\mathfrak t} (t))$ is the Weyl discriminant.
\end{itemize}
The complement of the strongly regular set is of $f'$-measure zero, and this completely describes $f'$. 

As is clear, the pushforward measure contains all the information about the stable orbital integrals of the test function $\Phi$, discussed in the lectures of Kaletha and Ta\"ibi \cite{Kaletha}, but provides a more economical way to talk about those, without having to fix Haar measures. Moreover, we will see that measures appear more naturally in the transfer operators beyond endoscopy; all formulas involving measures can be translated to formulas involving orbital integrals, but one has to remember to multiply with the absolute value of the Weyl discriminant, according to \eqref{pushfmeasure}. 

Nonetheless, there will also be instances when formulas are better expressed in terms of \emph{half-densities}, instead of measures. We will encounter this in our discussion of Hankel transforms (\S~\ref{sssHankel}), but also when we attempt to interpret some transfer operators in terms of quantization (\S~\ref{ss:transferrankone}).

Since there is no canonical notion of pushforward of a function or half-density, we need to invent one. This is what people typically do with orbital integrals, where they fix measures compatibly between the stable classes. In the setting above, we need to choose the Haar measures $dt$, as $T$ varies, in a compatible way. We do that by expressing such a Haar measure $dt$ as the absolute value of an invariant volume form $\omega_T$.\footnote{Fixing a measure $|dx|$ on $F$ and a volume form $\omega_X$ on a smooth variety $X$, we obtain a measure $|\omega_X|$ on $X(F)$, which in local coordinates looks $\omega_X = \phi(\underline{x}) dx_1 \cdots dx_n$ $\Rightarrow$ $|\omega_X| = |\phi(\underline{x})| |dx_1 |\cdots |dx_n|$. Of course, as is common, for the most part we will be writing $dx$ instead of $|dx|$, etc., when it is clear that it refers to a measure, not a differential form.} The set of such volume forms can be identified with the $F$-points of the top exterior power $\bigwedge^{\rk G} \mathfrak t_{\bar F}^*$. For any two maximal tori $T$, $T'$, if we fix an element over the algebraic closure that conjugates them, $g^{-1} T_{\bar F} g = T'_{\bar F}$, $g\in G(\bar F)$, this induces an isomorphism $\bigwedge^{\rk G} \mathfrak t_{\bar F}^* = \bigwedge^{\rk G} {\mathfrak t_{\bar F}'}^*$.  A compatible choice of Haar measures $dt$, $dt'$, then, is defined to be $dt = |\omega_T|$, $dt' = |\omega_{T'}|$, where $\omega_T$, $\omega_T'$ are $F$-rational invariant volume forms, and equal under the above $\bar F$-isomorphism of top exterior powers, up to an element of $\bar F^\times$ of absolute value $1$. 

Making such compatible choices allows us to talk about the stable orbital integrals of $\Phi$ as a function on the regular semisimple subset of $\Dfrac{G}{G}(F)$, up to a scalar choice that is independent on the orbit. We will have to live with this ambiguity, locally, and will denote by $p_! \Phi$ the corresponding ``pushforward function'' on the strongly regular semisimple set,  
\begin{equation}\label{pushffunction} p_!\Phi ([t]) = \SO_t (\Phi).
\end{equation}

Finally, for the half-density $\varphi = \Phi (dg)^\frac{1}{2}$, we will interpolate between \eqref{pushfmeasure} and \eqref{pushffunction}, and define, having fixed compatible choices of Haar measures, 
\begin{equation}\label{pushfdensity} p_!\varphi ([t]) = \SO_t (\Phi) |D(t)|^\frac{1}{2} (d[t])^\frac{1}{2}.
\end{equation}

\begin{remark}\label{remark-stable}
There is a Schwartz space of test measures $\mathcal S(\frac{G}{G})$ associated to the quotient (stack) of $G$ by conjugation, which corresponds to the zeroth homology of a Schwartz complex defined in \cite[\S~3]{SaStacks}. This is identified with the direct sum of coinvariant spaces $\mathcal S(G')_{G'}$ of the Schwartz spaces of all pure inner forms $G'$ of the group, where we abuse notation again to write $G$ for $G(F)$, etc. However, in the group case the reader does not need to worry about these definitions, since we are only interested in \emph{stable} coinvariants, and those are identified with the image of the pushforward map from the quasisplit form only. 

The fact that it is enough to consider one form of the group is a special feature of the group case.  For more general (say, smooth affine) $G$-varieties $X$, there is a well-defined pushforward map $\mathcal S(X/G)\to \operatorname{Meas}(X\sslash G)$, whose image will be called the space $\mathcal S(X/G)^\st$ of stable test measures for this quotient, but it cannot, in general, be identified with the coinvariant space of a single form. A simple example is Jacquet's relative trace formula for torus periods \cite{Jacquet-Waldspurger}, where $X=(T\backslash \PGL_2)^2$ for $T$ a nonsplit torus, and $G=\PGL_2$, acting diagonally; here, considering the quaternion division algebra is essential for obtaining the complete space $\mathcal S(X/G)$. Of course, one can always choose an embedding $G\hookrightarrow \GL_n$; then the space $\mathcal S(X/G)$ can be identified with the coinvariant space $\mathcal S(X\times^G \GL_n)_{\GL_n}$, and $\mathcal S(X/G)^\st$ is its pushforward to $X\sslash G = (X\times^G \GL_n)\sslash\GL_n$. 
\end{remark}

\subsection{Transfer operators}

In the (conjectural, in general) local Langlands correspondence, every tempered $L$-packet $\Pi_\phi$ for $G(F)$ gives rise to a stable character $\Theta_\phi$, defined as a suitable linear combination of the characters in the $L$-packet \cite{Kaletha}, and it is expected \cite[Conjecture 9.2]{Shahidi-Plancherel} that those stable characters span a (weak-*) dense subspace of the dual of $\mathcal S(\frac{G}{G})^\st$. Given a morphism between the $L$-groups of two quasisplit reductive groups, 
\[ \iota: {^LH} \to {^LG},\]
the local Langlands conjecture associates to every tempered $L$-packet of $H$ a tempered $L$-packet of $G$, hence to any tempered stable character $\Theta_\phi$ on $H$ a tempered stable character $\Theta_{\iota_*\phi}$ of $G$.

\begin{remark} \label{remarknontempered}
These statements should also hold without the ``tempered'' adjective; however, it may not be reasonable to expect the techniques of ``Beyond Endoscopy'' to capture character relations for nontempered $L$-packets; the reason is that it is the characters of standard modules, not of their Langlands quotients, which behave nicely in families. 
\end{remark}

The question that Langlands posed in \cite{Langlands-ST} is whether one can describe a ``transfer'' map between spaces of stable test measures
\begin{equation}\label{transfermap}
 \mathcal T_\iota: \mathcal S(\frac{G}{G})^\st \to \mathcal S(\frac{H}{H})^\st ,
\end{equation}
whose adjoint would give rise to the transfer of stable characters predicted by functoriality, i.e., 
\begin{equation}\label{transfercharacterization}
  \Theta_\phi \circ \mathcal T_\iota = \Theta_{\iota_*\phi}
\end{equation}
for every tempered $L$-packet $\phi$ of $H$.

We should clarify the assumptions and logical order here: One needs to be given the map $\Theta_\phi\mapsto \Theta_{\iota_*\phi}$ of stable local Langlands functoriality, and the density of stable characters, in order to make sense of \eqref{transfercharacterization} as a property characterizing the transfer map. This is indeed the case in the setting of \cite{Langlands-ST}, where Langlands focuses on the case $H=$ a torus and $G=\GL_2$, where local functoriality is explicitly known. However, in general, one could hope that the transfer map \eqref{transfermap} would be described first, as was the case for transfer factors in the theory of endoscopy, and then employed in a comparison of (stable) trace formulas to prove local Langlands functoriality. In any case, whether local functoriality is taken as input or not, the question of the existence and explicit calculation of the transfer \eqref{transfermap} is important for establishing global functoriality.

As in the case of endoscopy, the transfer \eqref{transfermap} is only part of the local input needed for a global comparison. One would additionally need to control the image of the transfer for some ``basic function'' (or ``basic measure'') at unramified places, and for its translates under Hecke operators. In endoscopy, this basic measure is the unit element of the Hecke algebra, and the statement about its transfer is the Fundamental Lemma. Beyond endoscopy, as we will see, one needs to work with more general ``basic measures.'' Indeed, according to Langlands' idea, one needs to insert $L$-functions into the trace formula in order to extract, via their poles, images of functorial lifts, and this is accomplished by modifying the basic measure. We will discuss this in the sections that follow. In any case, the fundamental lemma for the standard basic function (and the entire Hecke algebra) is a triviality if one knows the transfer relation \eqref{transfercharacterization}. In practice, however, as in the case of endoscopy, one needs to work in the reverse order: first establish the fundamental lemma for the Hecke algebra, then draw conclusions about character relations.

\subsection{Local transfer from tori to $\SL_2$}\label{ss:toriSL2}

The questions posed above have only been fully addressed in the case $H=T=$ a $1$-dimensional torus and $G=\SL_2$. The map of $L$-groups is 
\[\iota: {^LH} = \Gm\rtimes \Gal(E/F) \to \check G = \PGL_2,\]
where $E/F$ is the splitting field of the torus $T$, and we feel free to abbreviate $L$-groups by their relevant quotients in every case. The torus could be split, in which case $E=F$; otherwise, the map above can be identified with the projectivization of the orthogonal group embedded into $\PGL_2$.

In \cite{Langlands-ST}, Langlands addressed the question of transfer for this case. 
In the non-Archimedean case, he relied on the character formulas of Sally and Shalika (for local fields of characteristic zero and residual characteristic $\ne 2$). In this case, the stable characters of $L$-packets for $\SL_2$ are known, hence it makes sense to appeal to \eqref{transfercharacterization} for a characterization of the transfer map. The results of Langlands can be summarized as follows, both in the Archimedean and in the non-Archimedean case:

\begin{theorem}[Langlands]
 In the setting above, the transfer operator $\mathcal T_\iota$ satisfying \eqref{transfercharacterization} exists. 
\end{theorem}

This is proven in \S~2 of that article, working backwards from \eqref{transfercharacterization}, using explicit character formulas. As is common in the theory of the trace formula, Langlands works with orbital integrals of test functions, but the statements become simpler if we work with measures, instead. The formula \eqref{pushfmeasure} relating the pushforward of a measure $f = \Phi dg$ and stable orbital integrals, in the notation of \cite[(2.6)]{Langlands-ST} reads 
\[p_!f (a_G) = |\Delta(a_G)|\lambda(a_G) \operatorname{Orb}^\st(a_G, \Phi) da_G.\] 

Despite the case-by-case analysis used in the proof, the formula for the transfer operator is surprisingly uniform and simple. It is essentially a formula discovered many decades earlier by Gelfand, Graev, and Piatetski-Shapiro \cite[2.5.4.(7)]{GGPS}. The formula states that if the torus $T$ is nonsplit, and we denote by $\eta$ the quadratic character associated to its splitting field, then the correspondence $\chi \mapsto \Theta_\chi$ between characters of $T$ and stable characters for $\SL_2$ is given by 
\begin{equation}\label{GGPS}
 \Theta_\chi = \frac{2}{\Vol(T)} \frac{\eta(x)}{|x|} {\star_+} S_\chi.
\end{equation}
Here, $\star_+$ denotes additive convolution on the affine line, identified both with the Steinberg--Hitchin base of $\SL_2$ (through the trace map) and with the quotient $T\sslash \mu_2$ for $T$, where $\mu_2$ acts by inversion. The identification of the latter is by means of the trace composed with any embedding $T\hookrightarrow \SL_2$. Equivalently, fixing an isomorphism $T_{\bar F} = \Gm{_{\bar F}}$ over the algebraic closure, the coordinate $y$ on $\Gm$ corresponds to $x= y+ y^{-1}$ on $\mathbb A^1$. This gives a well-defined morphism $T\ni t\mapsto x(t) \in \mathbb A^1$ over $F$, and a well-defined ``discriminant'' function $D(t) = x(t)^2-4$, which over the algebraic closure is equal to $(y-y^{-1})^2$. Finally, $S_\chi$ is the pushforward to $\mathbb A^1$ of the measure $\chi(t) dt$ on $T$; the choice of Haar measure $dt$ does not matter, since it is cancelled by the factor $\Vol(T)$ in the denominator.

The formula \cite{GGPS} does not literally make sense when $T$ is split (the distribution $\frac{dx}{|x|}$ on the affine line is not well-defined, and $\Vol(T)$ is infinite), but the first factor in the convolution can be understood as a delta function at $0$, and then the formula remains true: $\Theta_\chi$ is equal to $S_\chi$ in this case, divided by a suitable Haar measure $dx$, by a well-known formula for the character of principal series representations; explicitly:
\[ \Theta_\chi(x) = \begin{cases} \frac{\chi(t) + \chi^{-1}(t)}{|t-t^{-1}|},& \mbox{ if } x= t+t^{-1}, \,\, t\in F^\times; \\ 0,& \mbox{otherwise.} \end{cases}\]

 The uniformity of the two cases was explained in \cite[\S 10.4]{SaTransfer2}, from which we can extract the following explicit formula for the transfer operator $\mathcal T_\iota$, applied to an element $f\in \mathcal S(\frac{G}{G})^\st$, written as $f(x) = \Phi(x) dx$, with $dx$ a measure on the real line: 

\begin{equation} \label{transfer-TSL2}
\mathcal T_\iota(f) (t) =  dt \cdot \lambda(\eta,\psi) \int_F \eta(uv)  \int_F \psi(v x(t) - uv) \Phi(u) \, du \,\, dv, 
\end{equation}
understood as Fourier transform of $\Phi$, followed by multiplication by the character $\eta$, followed by inverse Fourier transform. Here, $\psi$ denotes a fixed unitary character of the additive group $F$, and $\lambda(\eta,\psi)$ is a scalar ratio of gamma factors, see \cite[(7.14)]{SaTransfer2}. The Haar measure $dt$ on the torus is the one that can locally be written as $|D(t)|^{-\frac{1}{2}} dx(t)$, using the fact that the map $t\mapsto x(t)$ is a local homeomorphism away from the zeroes of the discriminant.

\subsection{Kuznetsov transfer from tori to $\SL_2$} \label{ss:Kuztori}

Let us now discuss a different type of transfer for this functorial lift, which uses the Kuznetsov formula for $\SL_2$ instead of the stable trace formula, and was initiated in the thesis of Venkatesh \cite{Venkatesh}. The local theory that follows was developed in \cite[\S~10]{SaTransfer2}. 

The functorial lift that we want to discuss belongs to the generalization of the Langlands conjectures afforded by the ``relative'' Langlands program. More specifically, in the discussion above we will replace the group $G=\SL_2$ by its Whittaker model, and adjoint quotient of the trace formula by the corresponding quotient for the \emph{Kuznetsov formula}. That is, instead of Schwartz measures on $G$ we will consider the space $\mathcal S((N,\psi)\backslash G)$ of Whittaker Schwartz measures, and instead of the stable coinvariant space $\mathcal S(\frac{G}{G})$ we consider the coinvariant space $\mathcal S((N,\psi)\backslash G/(N,\psi)) = \mathcal S((N,\psi)\backslash G)_{(N,\psi)}$. 

Here, $N$ is the upper triangular unipotent subgroup, identified with the additive group in the standard way, and $\psi:F\to \mathbb C^\times$ is a nontrivial unitary character. It will be convenient, at first, to think of the space $\mathcal S((N,\psi)\backslash G/(N,\psi))$ as a space of measures on the invariant-theoretic quotient $N\backslash G\sslash N$; there is no need to add the adjective ``stable'' here, as general orbits are stable, but because of the character $\psi$ we need to adopt some (non-canonical) conventions, in order to define a pushforward map 
\[\mathcal S((N,\psi)\backslash G/(N,\psi))\to  \operatorname{Meas}(N\backslash G\sslash N).\]
We will follow the convention of \cite[\S~2.2.1]{SaTransfer1}, according to which the image of a measure $\Phi dg \in \mathcal S(G)$ in $\mathcal S((N,\psi)\backslash G/(N,\psi))$ (where $\Phi$ is a Schwartz function on $G$) will be identified with the measure 
$\phi(\zeta) |\zeta| d\zeta$
on $\mathbb A^1\simeq N\backslash G\sslash N$, where $\phi$ is the function of orbital integrals
\[ \phi(\zeta)=\int_{F^2} \Phi \left(\begin{pmatrix} 1& x \\ &1\end{pmatrix} \begin{pmatrix} & -\zeta^{-1} \\ \zeta \end{pmatrix} \begin{pmatrix} 1& y \\&1 \end{pmatrix} \right) \psi^{-1}(x+y) dx dy.\]
The Haar measure $dg$ here is such that on the open Bruhat cell  it coincides with $|\zeta| d\zeta\, dx \, dy$, in the coordinates of the formula.

The dual group of the Whittaker space remains $\check G$, therefore the same map $\iota$ of $L$-groups as above predicts the existence of a transfer operator
\[ \mathcal T_\iota: \mathcal S((N,\psi)\backslash G/(N,\psi)) \to \mathcal S(T),\]
behaving well with respect to characters. For the Kuznetsov formula, ``characters'' refers to the ``relative characters,'' also called ``Bessel distributions/characters'' in this case. For a general discussion of relative characters, I point to the article of Beuzart--Plessis \cite{Beuzart-Plessis}. Recall that, unlike in the group case, the relative characters are only canonical up to a scalar depending on the representation, and a choice of Plancherel measure is needed to fix that scalar. For the Whittaker model we can use the group Plancherel measure, reducing the ambiguity to a scalar that does not depend on the representation. 

If we then denote by $J_\chi$ the Bessel character corresponding to the tempered $L$-packet corresponding to the unitary character $\chi$ of $T$ (actually, corresponding to the single element of the $L$-packet that is generic), the following theorem was proven in \cite[Theorem 5.0.1]{SaTransfer1}, \cite{SaTransfer2} (stated up to choices of Haar measures that we will not prescribe here):

\begin{theorem}
 There is a transfer map 
\[ \mathcal T_\iota: \mathcal S((N,\psi)\backslash G/(N,\psi)) \to \mathcal S(T)\]
 with the property that $T_\iota^*\chi = J_\chi$ for every unitary character $\chi$. Writing $f(\zeta) = \Phi(\zeta) d\zeta$, it is given by the formula 
  \begin{equation}\label{transfer-Venkatesh}
\frac{T_\iota(f) (t)}{dt} =  \lambda(\eta,\psi) \int \Phi\left(\frac{x(t)}{y}\right) \eta(x(t)y) |y|^{-1} \psi(y) dy.  
 \end{equation}
The function $T\ni t\mapsto x(t) \in \mathbb A^1$ and the scalar $\lambda(\eta,\psi)$ are the ones mentioned in \S~\ref{ss:toriSL2}, and $dt$ is a Haar measure on $T=T(F)$, (suitably) compatible with $d\zeta$. 
\end{theorem}

In other words, up to inverting the argument, the transfer is given by a simple Fourier transform of the measure $f\cdot \eta$. 

There is more to the above theorem than stated: The transfer map extends to a larger, ``nonstandard'' space of test measures 
\[\mathcal S_{L(\Ad,\eta |\bullet|)} (N,\psi\backslash \SL_2/N,\psi) \supset  \mathcal S (N,\psi\backslash \SL_2/N,\psi),\] 
which, at unramified non-Archimedean places, contains a distinguished ``basic vector,'' $f_{L(\Ad,\eta |\bullet|)}$, tailored to insert $L$-functions to the spectral side of the Kuznetsov formula. The insertion of $L$-functions is essential for the global comparison -- we will have more to say about it later.

Moreover, we have a fundamental lemma for this basic vector: its image under the transfer operator is equal, up to an explicit constant $c$, to the basic vector $1_{T(\mathfrak o)} \in \mathcal S(T)$. Furthermore, this fundamental lemma extends to the Hecke algebra: We have 
\begin{equation}\label{FL-torus}
 T_\iota(h\star f_{L(\Ad,\eta |\bullet|)}) = c\cdot \iota^*(h)\star 1_{T(\mathfrak o)},
\end{equation}
where $\iota^*$ is the map from the unramified Hecke algebra of $T$ to the unramified Hecke algebra of $G$ induced by the map of $L$-groups $\iota$ via the Satake isomorphism.

\subsection{Some early takeaways}

Some takeaways from our discussion of transfer operators for the lift from tori to $\SL_2$ are the following:

\begin{enumerate}
 \item At least in this low-rank case, the transfer operators are computable; most significantly, \emph{the formulas involve abelian Fourier transform}, despite the fact that they realize functorial lifts to the non-abelian group $\SL_2$. This is a very promising sign! It would be desirable, however, to achieve more conceptual understanding of these formulas, in order to be able to guess the transfer operators for higher-rank cases, where the functoriality conjecture is still open.
 \item The form of the transfer operators makes them suitable, at least in principle, for an application of the Poisson summation formula over the adeles. This was accomplished in Venkatesh's thesis \cite{Venkatesh} for the transfer operator \eqref{transfer-Venkatesh}, as we will see in \S~\ref{ss:Venkatesh}. Even if the operators have the form of Fourier transforms, proving a global Poisson summation formula in this setting is quite nontrivial, since the transforms are applied to nonstandard spaces of test functions/measures. We will discuss such global comparisons in Section \ref{sec:global}.
\item The formula \eqref{transfer-Venkatesh} for the transfer to the Kuznetsov formula is simpler than the one \eqref{transfer-TSL2} for the transfer to the stable trace formula. In fact, \eqref{transfer-TSL2} has been interpreted in \cite{SaTransfer2} as the composition of \eqref{transfer-Venkatesh} with another transfer operator, that we will encounter in Theorem \ref{Rudnick-local}. The Kuznetsov formula seems to be the nodal point for all transfer operators known to-date.
\end{enumerate}

\subsection{Comparisons involving the adjoint quotient in higher rank}

Generalizing the formula \eqref{GGPS} of Gelfand, Graev, and Piatetski-Shapiro to higher rank (in the non-Archimedean case) was the objective of the thesis of D.\ Johnstone \cite{Johnstone}. His work aims to study the stable transfer for a map of $L$-groups
\[{^LT}\to\PGL_n,\]
when $T$ is an $n$-dimensional elliptic torus. His work relies on character formulas and local character expansions for tame supercuspidal representations due to Adler, DeBacker, Murnaghan, and Spice, and is more complete when $n$ is a prime (different from the residue characteristic). The generalization of \eqref{GGPS} to higher rank is considerably complicated, and I point the reader to Johnstone's paper for the formula. An important question is whether the resulting transfer operator can be brought to a form where it could be shown to satisfy a Poisson summation formula on the Steinberg--Hitchin base, at least ``in principle.''

D.\ Johnstone and Z.\ Luo \cite{Johnstone-Luo} have studied transfer operators for another, very important family of functorial lifts, namely, the lifts corresponding to the symmetric power representations of the gual group of $\GL_2$:
\[ \Sym^n: \GL_2\to \GL_{n+1}.\]
In that paper, they compute a transfer operator or, rather, a distribution kernel, which pulls back ``most'' tempered characters of $\GL_2$ to the tempered characters of their $\Sym^n$-lift to $\GL_{n+1}$. The formula has a very simple form, but does not produce the desired outcome for twists of the Steinberg representation. It would be interesting to find the correct tweak of their formula that does that, and to rigorously show that it corresponds to a transfer operator between spaces of stable test measures.

\section{Same-rank comparisons and separation of the Arthur classes} \label{sec:samerank}

\subsection{Stable Kuznetsov--Selberg transfer for $\SL_2$} \label{ss:KuzSL2}

In his 1990 PhD thesis \cite{Rudnick}, Z.\ Rudnick showed how to use the Petersson formula -- that is, the simplified version of the Kuznetsov formula for holomorphic cusp forms -- to prove the Eichler--Selberg trace formula (the trace formula for holomorphic cusp forms).

In his thesis and subsequent work \cite{Altug1,Altug2,Altug3}, A.\ Altu\u{g} studied, among other questions, the problem of extracting the ``Ramanujan'' spectrum (that is, the complement of the 1-dimensional characters in the set of automorphic representations) from the Selberg trace formula for $\GL_2$. This was preceded by the work of Frenkel, Langlands, and Ng\^o \cite{FLN}, who sketched an argument for removing the contribution of the trivial representation, via an adelic Poisson summation formula.

All the aforementioned works are related, and have an underlying local theory. The local theory can be expressed in terms of a transfer operator
\[ \mathcal S((N,\psi)\backslash\SL_2/(N,\psi)) \to \mathcal S(\frac{\SL_2}{\SL_2})^\st\]
between spaces of test measures for the Kuznetsov and the stable trace formula for $\SL_2$. This transfer operator should give rise to a putative comparison between the stable trace formula and the Kuznetsov formula for $\SL_2$ (generalizing Rudnick's comparison), with an extra term corresponding to the trivial representation, since that is included in the spectrum of the former, but not of the latter. We will discuss this putative comparison, and its relation to the aforementioned works, in the global discussion in Section \S~\ref{ss:Altug}.

Consider the space of test measures for the Kuznetsov formula of the group $G=\SL_2$, $\mathcal S((N,\psi)\backslash G/(N,\psi))$, understood as a subspace of measures on the affine line with coordinate $\zeta$, as in  \S~\ref{ss:Kuztori}. Recall that $\mathcal S(\frac{G}{G})^\st$ denotes the space of stable test measures for the Selberg trace formula, i.e., the image of the pushforward $\mathcal S(G) \xrightarrow{\tr_*} \operatorname{Meas}(\mathbb A^1)$. The two spaces of measures on the affine line are completely different: The Kuznetsov space involves measures that are smooth of rapid decay away from a neighborhood of $\zeta=0$, and over non-Archimedean fields involves local Kloosterman sums (see \cite[\S~1.3]{SaRankone}). The Selberg space involves measures that are smooth of rapid decay outside of the discriminant locus $x=\pm 2$ (recall that the coordinate $x$, here, corresponds to the trace), and in neighborhoods of this locus behave like multiplicative characters, see \cite[\S 5.3]{SaRankone}. Nonetheless, one can describe a transfer operator between the two; to express it, denote, for any $s\in \mathbb C$, 
\begin{equation}\label{Ds}
D_s(x) = \psi(x) |x|^s d^\times x = \psi(x) |x|^{s-1} dx, 
\end{equation}
a measure on $F^\times$. 

\begin{theorem}[{\cite[Theorem 4.2.1]{SaTransfer1}}]\label{Rudnick-local}
 Consider the action of $F^\times$ on $F$, and denote by $\star$ the pushforward of measures under the action map (multiplicative convolution). Multiplicative convolution by $D_1$, $f\mapsto D_1\star f$, gives rise to an injection
 \[ \mathcal T: \mathcal S((N,\psi)\backslash G/(N,\psi)) \to \mathcal S(\frac{G}{G})^\st,\]
with the property that $\mathcal T^* \Theta_\Pi = J_\Pi$ for every tempered $L$-packet $\Pi$. Here, $\Theta_\Pi$ is the stable character of the $L$-packet (the sum of the characters over elements of the $L$-packet), and $J_\Pi$ is the Kuznetsov relative character (Bessel distribution) for the generic element of the $L$-packet (see \S~\ref{ss:Kuztori}). 
\end{theorem}

To explicate, this convolution operator applied to a measure of the form $f(\zeta)=\Phi(\zeta) d\zeta$ is the Fourier transform of the function $y\mapsto |y|^{-1} \Phi(y^{-1})$: 
\[ D_1\star f (x) = dx\cdot \int_y |y|^{-1} \Phi(y^{-1}) \psi(yx) dy.\]
The inverse operator sends a measure $f'(x)=\Phi'(x) dx$ in the trace variable $x$ to the measure 
\begin{equation}\label{inversetransfer} |\zeta|^{-1} d\zeta\cdot \int_x \Phi'(x) \psi^{-1}(\frac{x}{\zeta}) dx
\end{equation}
in the variable $\zeta$ for the Kuznetsov quotient. 

This additive Fourier transform of $\Phi'$ is prominent in a big part of the literature of ``Beyond Endoscopy'' \cite{FLN, Altug1, Altug2, Altug3}, where the Kuznetsov formula is not explicitly invoked. One of the goals of these papers has been the removal of the trivial representation from the stable trace formula of $\SL_2$.\footnote{Most of these papers concern $\GL_2$ or $\PGL_2$, but for simplicity we will restrict ourselves to $\SL_2$; the passage to $\GL_2$ is rather innocuous.} The relevant local observation here is that the Fourier transform of the measure $f' \in \mathcal S(\frac{G}{G})^\st$, evaluated at $0$ (that is, at $\zeta^{-1}=0$, in the coordinates above), is simply the integral of $f'$, i.e., the trace of $\tilde f'$ at the trivial representation (where $\tilde f' \in \mathcal S(G)$) is any preimage of $f'$. How to utilize this easy observation globally is one of the questions that the aforementioned papers tried to address; it will be discussed in \S~\ref{ss:Altug}.

The theorem above reveals more about the nature of this Fourier transform: It is the operator of functoriality from the Selberg to the Kuznetsov formula. Note that \eqref{inversetransfer} does not preserve rapid decay in the variable $\zeta$: for example, it makes sense to evaluate the integral appearing in \eqref{inversetransfer} at $\zeta=\infty$; as mentioned, this computes the total mass of the measure $f'$ (the character of the trivial representation). This means that, to turn the transfer operator $\mathcal T$ to an isomorphism, we need to apply it to an \emph{enlargement}
\[\mathcal S_{L(\Ad,1)}((N,\psi)\backslash G/(N,\psi))\supset \mathcal S((N,\psi)\backslash G/(N,\psi))\]
of the standard space of test measures, which has both a representation-theoretic and an arithmetic role. Both can be understood more easily in terms of global trace formulas: A direct comparison between the stable Selberg and Kuznetsov trace formulas for $\SL_2$ would be impossible, because:
\begin{enumerate}
 \item The former contains a contribution from the trivial representation, which does not appear in the latter.
 \item The traces of Hecke operators on the latter are weighted by factors of the form 
 \[ \frac{(\mbox{\small{local factors at }}S)}{L^S(\pi, \Ad, 1)},\]
 where $S$ is a finite number of places, and $L^S(\pi, \Ad, 1)$ denotes the partial adjoint $L$-function of the automorphic (say, cuspidal) representation $\pi$. 
\end{enumerate}

Using the enlarged space $\mathcal S_{L(\Ad,1)}((N,\psi)\backslash G/(N,\psi))$ of test measures resolves the two issues, at least in principle. We will not get into the explicit description of this space here
 (see \cite[\S~2.2]{SaTransfer1}), because it is not very conceptual, and certainly not well-understood in higher rank. What is understood, and is the source of the notation, is that at non-Archimedean places it contains a ``basic vector'' $f_{L(\Ad,1)}$ which is a generating series for the local unramified adjoint $L$-value at $1$. More precisely, assume that $F$ is non-Archimedean with residual degree $q$, and consider the following formal series in the representation ring of $\check G=\PGL_2$:
 \[\hat f:=\sum_{n\ge 0} q^{-n} \Sym^n \Ad.\]
 The trace of a semisimple conjugacy class $c$ in $\check G$ on $\hat f$ computes the local adjoint $L$-value $L(\pi,\Ad,1)$ for the irreducible unramified representation $\pi$ of $G=G(F)$ with Satake parameter $c$. The inverse Satake transform $f$ of $\hat f$ is a series of elements in the Hecke algebra of $G$, which converges as a measure on $G$. The corresponding series of pushforwards to $\mathcal S((N,\psi)\backslash G/(N,\psi))$ converges as a measure on the affine line, to the measure that we call the ``basic vector'' of $f_{L(\Ad,1)}\in \mathcal S_{L(\Ad,1)}((N,\psi)\backslash G/(N,\psi))$. 
 
Theorem \ref{Rudnick-local} extends to identify (up to an explicit local zeta factor) the transfer of this basic vector to the basic vector of $\mathcal S(\frac{G}{G})$ (the image of the unit  in the Hecke algebra), see \cite[Theorem 4.2.1]{SaTransfer1}. 

\subsection{Stable Kuznetsov--Arthur transfer for $\SL_n$} \label{ss:KuzSLn}

We briefly discuss possible generalizations to $G=\SL_n$, for larger values of $n$. 

The question originally posed in this setting is how to remove the nontempered contributions from the stable trace formula of $G$, but one could also formulate the deeper question of how to compare the stable Arthur and Kuznetsov trace formulas. In \cite{Arthur-stratification}, Arthur formulates a hypothesis, that the non-tempered contributions to the trace formula can be obtained by evaluating the Fourier transform of a test measure $f\in \mathcal S(\frac{G}{G})^\st$ at certain subvarieties of the ``dual'' to the Steinberg--Hitchin base $\Dfrac{G}{G} \simeq \mathbb A^{n-1}$. 

In ongoing work, the author with Chen Wan are developing a full comparison between the Kuznetsov and stable Arthur--Selberg trace formulas. This comparison uses local transfer operators described by integral formulas which have been explicitly evaluated in low rank, and paint a different picture than the one obtained by treating $\Dfrac{G}{G}$ as an $(n-1)$-dimensional vector space. 

For simplicity, I only describe the case of $n=3$: Consider the additive character $\psi: F\to \mathbb C^\times$ as a generic character of the upper triangular unipotent $N\subset \SL_3$ in the standard way (i.e., applying it to the sum of the entries right above the diagonal), and use the matrices
\[\begin{pmatrix} && \zeta_3 \\ & \zeta_2 \\ \zeta_1 \end{pmatrix}\]
(with $\zeta_1\zeta_2\zeta_3=1$)
to represent the space of test measures $\mathcal S((N,\psi)\backslash G/(N,\psi))$ for the Kuznetsov formula as a space of measures on a $2$-dimensional space, as in \S~\ref{ss:Kuztori}. We identify the matrix above with the element $\diag(\zeta_1, \zeta_2,\zeta_3)$ of the Cartan group $T = B/N$, taking $B$ to be the Borel subgroup of upper triangular matrices; we denote the simple roots by $\alpha_1= \frac{\zeta_1}{\zeta_2}$ and $\alpha_2 = \frac{\zeta_2}{\zeta_3}$.

For the Steinberg--Hitchin base $\Dfrac{\SL_3}{\SL_3}$ we will use \emph{rational} coordinates $(u_1, u_2, u_3)$, with $u_1(g) = \tr g$, $u_2(g) = \frac{\tr \wedge^2 g}{\tr t}$, $u_3 = \frac{\det(g) }{\tr \wedge^2 g} $. (Obviously, $\det(g)=1$ and $u_1 u_2 u_3 =1$, but the formulas that follow carry over verbatim to $\GL_3$.) In upcoming work with C.~Wan, we describe a transfer operator 
 \[ \mathcal T: \mathcal S((N,\psi)\backslash G/(N,\psi)) \to \mathcal S(\frac{G}{G})^\st,\]
with the property that $\mathcal T^* \Theta_\Pi = J_\Pi$ for every tempered $L$-packet $\Pi$ (notation as in Theorem \ref{Rudnick-local}), and given by the following formula:
\begin{equation}\label{SWaneq}
 \mathcal Tf = D'_{\check\alpha_1, 1} \star D_{\check\alpha_2,1} \left(\psi^{-1}(e^{\alpha_1}) \psi(e^{\alpha_2}) D_{\check\alpha_1+\check\alpha_2,2} \star f \right).
\end{equation}
Here, for $s\in \mathbb C$ and a cocharacter $\lambda:\Gm\to T$, we denote by $D_{\lambda, s}$ the pushforward of the measure $D_s$ of \eqref{Ds} to $T$; we denote by $D_{\lambda, s}'$ the pushforward of the analogous measure, when $\psi$ is replaced by $\psi^{-1}$; and we denote by $\psi(e^\alpha)$ the operator of multiplication by $\psi$ composed with a character $\alpha: T\to \Gm$.\footnote{The exponential notation is used because the character group is written additively.}

The form of the transfer operator (including the choice of coordinates) may seem, and is, mysterious, but at least some degenerate limit of it, when we let the spaces to degenerate to their asymptotic cones, can be understood along the lines of \cite[\S~4.3, \S~5]{SaTransfer1}; in that limit, the factor $\psi^{-1}(e^{\alpha_1}) \psi(e^{\alpha_2})$ would disappear, and we just have Fourier convolutions along the positive coroot cocharacters. In any case, an important observation is that it is not a simple, $2$-dimensional Fourier transform on the $2$-dimensional Steinberg--Hitchin base that produces the transfer operator to the Kuznetsov formula (and hence, the separation of the non-tempered spectrum, which will again be responsible for a nonstandard space of Kuznetsov test measures, if we invert the operator $\mathcal T$), but a $3$-dimensional Fourier transform (together with the little-understood ``correction factor'' $\psi^{-1}(e^{\alpha_1}) \psi(e^{\alpha_2}$), along directions determined by the positive coroots. Higher-rank calculations for $\GL_n$ also support the idea that the relevant transfer operator involves $\frac{n(n-1)}{2}$ Fourier convolutions, together with correction factors.

\subsection{Transfer between relative trace formulas in rank 1, and a symplectic interpretation} \label{ss:transferrankone}

The comparison of Theorem \ref{Rudnick-local} generalizes to arbitrary affine, homogeneous spherical varieties $X$ whose dual group is $\SL_2$ or $\PGL_2$, up to possibly replacing the variety with a finite cover. That is (up to such replacement), there is a transfer operator 
\[ \mathcal T: \mathcal S((N,\psi)\backslash G^*/(N,\psi))  \to \mathcal S((X\times X)/G)^\st, \]
from Kuznetsov test measures on the group $G^*=\PGL_2$ or $\SL_2$ (depending on the dual group of $X$), and stable test measures for the relative trace formula of $X$. The case of the Selberg trace formula (Theorem \ref{Rudnick-local}) corresponds to the special case $X=\SO_3\backslash\SO_4$, but it is quite remarkable that essentially the same transfer operator, with small modifications, like replacing the distribution $D_1$ by $D_s$ for an appropriate $s$ (which, in turn, has interesting connections to the $L$-value associated with the spherical variety $X$), gives rise to the pertinent operators of functoriality. 

We will recall the theorem of \cite{SaRankone} for the family of spaces $X=\SO_{2n-1}\backslash\SO_{2n}$ only, where the transfer operator is given by an $1$-dimensional Fourier transform, but several more cases are treated in loc.cit., with a $2$-dimensional Fourier transform. The benefit of discussing this family is that the result was subsequently reproven by W.\ T.\ Gan and X.\ Wan \cite{Gan-Wan} using the theta correspondence, and their work allows us to formulate statements about characters (while \cite{SaRankone} only established the geometric transfer). 

For the theorem that follows, we can take $\mathcal S((X\times X)/G)^\st$ to be the pushforward of Schwartz measures from $X\times X$ to the affine line, where $X$ is the unit ``sphere'' on a \emph{split} quadratic space of dimension $2n$, and the affine line is identified with the quotient $(X\times X)\sslash G$ via the quadratic pairing. (As in the case of $\SL_2$, this is enough, and one does not obtain a larger space of test measures by considering pure inner forms.) We will formulate functoriality for characters in terms of the theta correspondence, following \cite{Gan-Wan}; for this, we need to choose the data of the theta correspondence so that the theta lift of an $(N,\psi)$-generic tempered representation on $\SL_2$ to $\SO_{2n}$ is distinguished by the split form of $\SO_{2n-1}$; see \cite[Proposition 5.1]{Gan-Wan} for details.

\begin{theorem}[S., Gan--Wan] \label{rankone}
Let $G^*=\SL_2$, $G=\SO_{2n}$, and $X=\SO_{2n-1}\backslash \SO_{2n}$, as above. The operator 
\[ \mathcal T f (x) = |x|^{n-2} D_{3-n}\star f (x)\] 
gives rise to an injection
 \[ \mathcal T: \mathcal S((N,\psi)\backslash G^*/(N,\psi)) \to \mathcal S((X\times X)/G)^\st,\]
with the property that $\mathcal T^* I_{\theta(\Pi)} = J_\Pi$ for every tempered $L$-packet $\Pi$ of $\SL_2$. Here, $J_\Pi$ is the Bessel character as in Theorem \ref{Rudnick-local}, and  $I_{\theta(\Pi)}$ is a relative character for the theta-lift of the $(N,\psi)$-generic element of $\Pi$ to $G$.
\end{theorem}

For the normalization of the relative character $I_{\theta(\Pi)}$, I point the reader to \cite[Theorem 8.1]{Gan-Wan}; suffice it to say here that it is a normalization compatible with the Plancherel formula for $L^2(X)$. This theorem also has an extension, describing an enlargement of the standard space of Kuznetsov test measures, such that $\mathcal T$ becomes an isomorphism
 \begin{equation}\label{transferisomorphism} \mathcal T: \mathcal S_{L(\Ad,n-1)}((N,\psi)\backslash G^*/(N,\psi)) \xrightarrow\sim \mathcal S((X\times X)/G)^\st.
 \end{equation}
 The inverse of this isomorphism takes the basic vector on the right hand side (the image of the probability measure which is equal to the characteristic function on $X\times X(\mathfrak o)$ times an invariant measure) to an explicit multiple of the generating series of the local unramified $L$-value $L(\Ad,n-1)$.

The uniform nature of these transfer operators is enticing, but mysterious. What are we to make of the fact that functorial transfer for homogeneous spaces of non-abelian reductive groups are given by Fourier convolutions -- operators of abelian nature -- on spaces of (stable) test measures? An interpretation using symplectic geometry and the idea of geometric quantization was proposed in \cite{SaICM}: According to this idea, the ``cotangent stacks'' associated to the quotients $\mathfrak Y=(N,\psi)\backslash G^*/(N,\psi)$ and $\mathfrak X = (X\times X)/G$ are, in some sense, ``birational'' to each other. More precisely, the Kuznetsov cotangent stack $T^*\mathfrak Y$ can be canonically identified with the group scheme 
\[ J_X = \left(\Res_{\mathfrak a_X^*/\mathfrak c_X^*} T^* A_X\right)^{W_X}\]
of regular centralizers in $\SL_2$; here,  $A_X=$ the (universal) Cartan of $\SL_2$, $W_X=$ its Weyl group, and $\mathfrak c_X^* = \mathfrak a_X^*\sslash W_X$, which is isomorphic to the affine line. I point the reader to \cite{SaICM} for more details and definitions. On the other hand, using work of Knop \cite{KnWeyl, KnAut}, one can show that there is a canonical diagram
\[ J_X \leftarrow J_X\times_{\mathfrak c_X^*} T^*X/G \to T^*\mathfrak X,\]
where the left arrow is generically an isomorphism, and the right arrow has dense image. One can interpret this by saying that the cotangent stack of $\mathfrak X$ is ``birational to, but a bit larger than'' the Kuznetsov cotangent stack $T^*\mathfrak Y$. Morally, Theorem \ref{rankone} and its extension \eqref{transferisomorphism} are the quantum versions of this fact. In \cite{SaICM}, an interpretation of the transfer operators is given, as ``canonical intertwiners'' between two different ``geometric quantizations'' of these stacks. 

This explanation is purely phenomenological at this point, but it would be nice to adapt it to all other available examples discussed in this survey; in particular, it would be nice to have an explanation based on symplectic geometry of the mysterious appearance of Fourier transforms in problems of functorial transfer.

\section{Schwartz spaces encoding $L$-functions} \label{sec:Schwartz}

In the introduction, we discussed the necessity of introducing ``nonstandard test functions'' (or measures) into the trace formula, in order to weigh the spectral side by $L$-functions. We already encountered nonstandard spaces of test measures in the local theory of the preceding sections; for example, for the full ``Beyond Endoscopy'' matching between the Kuznetsov and the stable trace formula of $\SL_2$ \eqref{transferisomorphism}, we had to use a nonstandard space $\mathcal S_{L(\Ad,s)}((N,\psi)\backslash G^*/(N,\psi))$ of test measures for the Kuznetsov formula. It is time to take a closer look at such spaces.

A basic example of nonstandard test measures is the space of Schwartz measures on the space $M_n$ of $n\times n$-matrices (over a local field $F$), considered as measures on $G=\GL_n$. The local theory of Godement and Jacquet \cite{GJ} shows that this space ``encodes'' the standard $L$-function; at non-Archimedean $F$ with ring of integers $\mathfrak o$ and residual degree $q$, this means: 
\begin{enumerate}
 \item the integral against the matrix coefficient of any irreducible smooth representation $\pi$, twisted by powers of the determinant:
 \begin{equation}\label{GJintegral}
  \mathcal S(M_n) \otimes \tilde\pi \otimes \pi \otimes |\det|^s \ni f \otimes \tilde v \otimes v \mapsto  \int f(g) \langle \tilde v, \pi(g) v\rangle |\det(g)|^s 
 \end{equation}
 is a rational function in $q^{\pm s}$, and has image equal to the fractional ideal generated by the standard (local) $L$-function, $L(\pi, \Std, s+\frac{n+1}{2}) \CC[q^{\pm s}]$;
 \item there is an unramified ``basic vector'' $f_\Std \in \mathcal S(M_n)^{G(\mathfrak o)\times G(\mathfrak o)}$ such that, for $\pi$ unramified, the integral \eqref{GJintegral}, against the normalized unramified matrix coefficient of $\pi$ is equal to $L(s+\frac{n+1}{2} , \pi, \Std)$. Up to the choice of additive Haar measure on $M_n$, the basic vector, here, is the characteristic function of $M_n(\mathfrak o)$. 
\end{enumerate}

Moreover, the functional equation of this $L$-function follows from a local and a global ingredient: the Fourier transform on $\mathcal S(M_n)$, and the Poisson summation formula.

It was an idea due to Braverman and Kazhdan \cite{BK2} that one should be able to define more general ``Schwartz spaces'' $\mathcal S_r$ of measures -- or, more conveniently, to center everything at $s=\frac{1}{2}$, spaces $\mathcal D_r$ of half-densities -- on the $F$-points of a reductive group $G$, possessing the properties above when we replace the standard representation of $\GL_n$ by an arbitrary algebraic representation $r$ of the $L$-group ${^LG}$.  They proposed that, at least when $r$ is irreducible, the space $\mathcal D_r$ should be obtained as the sum of the standard Schwartz space $\mathcal D_1=\mathcal D(G)$ and its image under a conjectural ``$r$-Fourier tranform,'' an $G\times G$-equivariant isometry on $L^2(G)$. We will discuss those conjectural Fourier transforms in the next section. 

In particular, when $G$ is unramified over a non-Archimedean field $F$, the space $\mathcal D_r$ should contain the following ``basic vector'' $f_r$. For technical reasons, let us assume, as in the Godement--Jacquet case, that there is a central cocharacter $\Gm\to {^LG}$ which, composed with $r$, gives the standard scaling action on $V$. Dually, this translates to a character $G\to \Gm$ which we will call ``determinant.'' 
If $\tilde h_n$ is the element of the unramified Hecke algebra of $G$ that corresponds to the representation $\Sym^n r$ under the Satake isomorphism, then, up to a suitable Haar half-density on $G$, the basic vector of $\mathcal D_r$ should be
\[f_r = \sum_{n\ge 0} q^{-\frac{n}{2}} h_n.\]
Finally, Braverman and Kazhdan observed that the support of this half-density  would have compact closure in a certain affine (but possibly singular) $G\times G$-equivariant embedding $G\hookrightarrow G_r$. Such embeddings automatically extend the multiplication operation on the group, turning them into monoids.

Their ideas were picked up by others about a decade later, and gave rise to several developments in the field. Motivated by his work with Frenkel and Langlands \cite{FLN,FN} on ``Beyond Endoscopy,'' and the need to insert $L$-functions into the trace formula, Ng\^o \cite{Ngo-PS} conjectured that, over a local field in equal characteristic, $F=\FF_q((t))$, the basic function\footnote{We feel free to switch between half-densities and functions or measures by dividing/multiplying by a suitable eigen-half-density. There is something slightly subtle here, in that this eigen-half-density is not necessarily invariant. For example, in the Godement--Jacquet case, the space $\mathcal D_{\Std}$ should be defined as the product of the space of Schwartz functions by $(dx)^\frac{1}{2}$, where $dx$ is \emph{additive} Haar measure on $M_n$. In other cases, with $r$ irreducible, the space $\mathcal D_r(G)$ should contain the ``basic half-density'' $\Phi_r(g) |\det g|^{\frac{1}{2}+s_r}$, where $\Phi_r$ is the ``IC function'' of the $L$-monoid $G_r$ determined by the heighest weight $\lambda_r$ of $r$ (see \cite{BNS}), and $s_r =\left< \rho_G, \lambda_r\right>$, where $\rho_G$ is the half-sum of positive roots of $G$. (The reason for this normalization lies in the calculation of the Godement--Jacquet integral of $\Phi_r$ -- see \cite{BNS-erratum}.)} $f_r$ should arise as a Frobenius trace from the intersection complex of the arc space $L^+G_r$, for $G_r$ now defined over $\FF_q$ (so that $L^+G_r(\FF_q)=G_r(\mathfrak o)$). The idea that basic functions representing $L$-functions should arise as ``IC functions'' associated to possibly singular affine spherical varieties had previously appeared in \cite{SaRS}, and was inspired by some other  work of Braverman and Kazhdan, on normalized Eisenstein series \cite{BK1,BK3}, which has its roots in the theory of geometric Eisenstein series and ideas of Drinfeld \cite{BG-Eisenstein}. In any case, a proper definition of the IC function was not available until \cite{BNS}, where Ng\^o's conjecture was proven.

A complete such Schwartz space $\mathcal S_r$ should also represent the $L$-function associated to $r$ for ramified representations, as well, as in the Godement--Jacquet case. A recent preprint of Bezrukavnikov, Braverman, Finkelberg, and Kazhdan \cite{BBFK} constructs the Iwahori-invariants of such a Schwartz space, proving the desired relation to $L$-functions.

Assuming, now, the existence of $r$-Schwartz spaces of functions, half-densities, or measures on our reductive group, we can consider various pushforwards. For example, unipotent integration against a generic character $\psi^{-1}$ of the maximal unipotent group $N$ gives rise to a morphism to the space of Whittaker functions or measures,
\[ \mathcal S_r(G) \to C^\infty((N,\psi)\backslash G),\]
whose image we can denote by $\mathcal S_r((N,\psi)\backslash G)$. It is a ``nonstandard Schwartz space of Whittaker functions,'' whose basic vector (the image of the vector $f_r$ above) already appeared in the musings of Piatetski--Shapiro around Poincar\'e series in the 80s \cite{PS-invariant}. Taking a two-sided pushforward, we obtain a nonstandard space $\mathcal S_r((N,\psi)\backslash G/(N,\psi))$ of test vectors for the Kuznetsov formula. This space \emph{should} coincide with the spaces that we encountered in Sections \ref{sec:localtori} and \ref{sec:samerank}, in particular cases.

\medskip

The aforementioned ``Schwartz'' spaces of measures could, more generally, be understood in the context of the ``relative'' Langlands program, introduced in the article of Beuzart-Plessis in these proceedings \cite{Beuzart-Plessis}. Indeed, the affine monoids that we encountered above are special cases of \emph{affine spherical $G$-varieties} $X$, and those should also be equipped with Schwartz spaces, often related to some $L$-function. To explain how, we need to make an adjustment to our understanding of the group case: First of all, we will now denote the group of interest by $X=H$, with $G=H\times H$ acting on it. Secondly, to fit into the general setting, the putative Schwartz space of the monoid $H_r$ should be understood as being related, formally, to $L(\pi, r) \sqrt{L(\pi, \Ad, 1)}$; for example, $\mathcal S(H)$ is related to $\sqrt{L(\pi,\Ad, 1)}$, not the trivial $L$-value. To give a meaning to this relationship, one needs to think of the Plancherel decomposition of the basic vector, which in the group case takes the form (for a suitable choice of Haar measure)
\begin{equation}\label{Plancherel-group} \Vert 1_{H(\mathfrak o)}\Vert^2 =  \int_{\check A_H^0/|W_H|} L(\chi,\Ad,1) d_{\text{Weyl}}\chi.
\end{equation}
Here, $\check A_H^0/|W_H|$ is the set of conjugacy classes in the compact real form of the dual group of $H$, and $d_{\text{Weyl}}\chi$ is the pushforward of the probability Haar measure on this compact group. This, of course, is just the well-known formula for the unramified Plancherel measure, due to Macdonald \cite{Macdonald}, and, in more generality, to Harish-Chandra \cite{Waldspurger-Plancherel}. 

The appearance of the adjoint (local) $L$-value in the Plancherel formula above is of significance in various aspects of the Langlands program. For example, it is related to the fact that the corresponding global $L$-value is equal, essentially, to the Petersson square norm of a normalized newform. ``Normalized,'' here, means that the first Fourier/Whittaker coefficient is $1$; the local counterpart of this this ratio between $L^2$-norms and Whittaker coefficients is the Plancherel density of the basic vector $W_0 \in \mathcal S((N,\psi)\backslash G)$ of the Whittaker model of $G$, which, in contrast to \eqref{Plancherel-group}, reads 
\begin{equation}\label{Plancherel-Whittaker} \Vert W_0\Vert^2=  \int_{\check A_G^0/|W_G|}  d_{\text{Weyl}}\chi.
\end{equation}
In the geometric Langlands program, the appearance of the adjoint $L$-value in \eqref{Plancherel-group} is understood to be related to the \emph{derived geometric Satake equivalence} \cite{BezFin}.

As is explained in the paper of Beuzart-Plessis in these proceedings, the Schwartz spaces of other spherical $G$-varieties $X$ admit similar Plancherel decompositions, with the Plancherel density of the ``basic vector'' encoding a different $L$-value, as in \eqref{Plancherel-group}, \eqref{Plancherel-Whittaker}. Moreover, the generalized form of the Ichino--Ikeda conjecture, proposed in \cite{SV}, relates local Plancherel densities to global period integrals of automorphic forms -- the same integrals that appear on the spectral side of the relative trace formula. Therefore, the insertion of $L$-functions into the Arthur--Selberg trace formula can be viewed as a relative trace formula for the monoid $H_r$. From this point of view, considering the Beyond Endoscopy program in the more general setting of the relative trace formula appears most natural.

\medskip

The Schwartz spaces of spherical varieties, when those are smooth, are well-defined. To complete the discussion of this section, let us return to the idea that singular spherical (and, maybe, more general) varieties also possess meaningful Schwartz spaces, which we already encountered in the case of monoids. We mention two classes of spaces where the basic vectors of these putative Schwartz spaces have been studied for $\mathfrak o = \mathbb F_q[[t]]$, and related to $L$-functions:

The first is the case of spaces of the form $X = \overline{U\backslash G}$ (or $X=\overline{[P,P]\backslash G}$), where $P\supset U$ denote a parabolic subgroup with its unipotent radical, and the line denotes ``affine closure.'' The basic functions of those spaces were studied in \cite{BFGM}. We will discuss them further in \S~\ref{ss:Fourier-nonlinear}, to highlight the nature (and absence!) of the conjectural nonlinear Fourier transforms.

Finally, \cite{SaWang} studied the basic vectors of Schwartz spaces of spherical varieties $X$ with dual group $\check G_X=\check G$. In all of those cases, the basic function, defined as the Frobenius trace on the intersection complex of the arc space of $X$ (or rather, on finite-dimensional formal models of that), admits a Plancherel decomposition as in \eqref{Plancherel-group}, \eqref{Plancherel-Whittaker}, with density equal to some local $L$-value that depends on the geometry of $X$.

\section{Fourier and Hankel transforms} \label{sec:Hankel}

\subsection{Desiderata} 

Let $X=\Mat_n$ be the space of $n\times n$-matrices, under the action of $\GL_n^2/\Gm$ by left or right multiplication (which we will write as a right action). Let $X^\vee$ be the dual vector space. If we fix a local field $F$ and a unitary additive character $\psi: F\to \mathbb C^\times$, we obtain a Fourier transform, an equivariant isomorphism
\[ \mathcal F: \mathcal D(X) \xrightarrow\sim \mathcal D(X^\vee)\]
between spaces of Schwartz half-densities on $X(F)$ and $X^\vee(F)$, which extends to an $L^2$-isometry. Given that both spaces contain $G=\GL_n$ as a canonical open subset, we can think of these half-densities as living on $G=G(F)$, and then Fourier transform is given by (multiplicative) convolution with the measure 
\[ \gimel_\Std(g) = \psi(\tr(g)) |\det g|^\frac{n}{2} dg,\]
where $dg$ is the Haar measure on $G$, depending on the choice of $\psi$, that makes the transform an $L^2$-isometry.

The same conjugation-invariant measure can be used to convolve functions and measures on $G$, and by the theory of Godement and Jacquet we have (essentially, by definition)
\begin{equation}\label{gamma-standard} \gimel_\Std \star \Theta_\pi = \gamma(\frac{1}{2},\pi, \Std, \psi) \Theta_\pi,
\end{equation}
for $\Theta_\pi$ the character (or, for that matter, any matrix coefficient) of an irreducible representation $\pi$ of $G$, and with $\gamma(\pi,r,\psi)$ the gamma factor of the functional equation,
\[ \gamma(s,\pi, r, \psi) = \frac{\epsilon(s,\pi, r,\psi) L(s,\tilde\pi, r)}{L(s,\pi, r) }.\] Here, we need some regularization to make sense of the convolution above (namely, using the Godement--Jacquet zeta integrals to meromorphically extend characters into the duals of Schwartz spaces of measures on $X$ and $X^\vee$), and the relation makes sense for almost all $\pi$, as they vary in families of the form $\pi\otimes |\det|^s$, $s\in\CC$.

Now, fix a reductive group $G$, for simplicity split, and an irreducible representation $r: \check G\to\GL(V)$ of its Langlands dual group. Braverman and Kazhdan \cite{BK2,BK4} considered the question of whether it is possible to describe an invariant distribution $\gimel_r$ on $G$, satisfying the analogous property
\begin{equation}\label{gamma-r} \gimel_r \star \Theta_\pi = \gamma(\frac{1}{2}, \pi,r,\psi) \Theta_\pi.
\end{equation}

Ideally, as discussed in Section \ref{sec:Schwartz}, there would be entire ``$r$-Schwartz spaces'' of half-densities $\mathcal D_r(G)$, with convolution by $\gimel_r$ defining an $L^2$-isometric isomorphism, the ``$r$-Fourier transform,''
\[ \mathcal F_r : \mathcal D_r(G) \xrightarrow\sim \mathcal D_{r^\vee}(G),\]
where $r^\vee$ is the dual representation. 

\subsection{The case of finite fields}

Although our goal is to have a definition of the $r$-Fourier transform over local fields (and, eventually, some type of Poisson summation formula over the adeles of a global field), Braverman--Kazhdan \cite{BK2,BK4} first focused on the (rather nontrivial) ``toy case,'' when $G$ is defined, and split, over a finite field. The study of this case was later advanced and completed by Cheng--Ng\^o and Chen \cite{CN,Chen}. We will only highlight the main points of this theory, pointing the reader to the original references, as well as \cite{Ngo-Takagi} for more details.

Let $G$ be a connected, split reductive group over a finite field $\FF=\FF_q$. The basic premise is that the distribution $\gimel_r$ on $G(\FF)$ should be obtained as the Frobenius trace of an irreducible perverse sheaf $\mathcal J_r$ on $G$ (up to a cohomological shift, which we will ignore below). Moreover, Braverman and Kazhdan constructed a candidate for this sheaf; the starting point is the analogous sheaf for the abelian case, namely, the case of the Cartan group $T$ of $G$, when we restrict the representation $r$ to its dual torus $\check T\subset \check G$; we will describe the construction below.

Of course, one also needs to define the analogs of gamma factors $\gamma(\pi,r)$ for the case of finite fields. This is done in \cite[\S~9]{BK2}, first for $\GL_n$, using the analog of the Godement--Jacquet construction, and finite-field Fourier transforms. (Here, one can use $\overline{\QQ_l}$-coefficients, for compatibility with the sheaf-to-function constructions that follow.) The gamma-factors are then transferred to any (split reductive) group $G$ via a version of ``functoriality for Deligne--Lusztig packets'' (= the set of irreducible constituents of a Deligne--Lusztig representation), via the map of dual groups $r: \check G \to \GL_n$. I point to the original reference for the definition, mentioning only that it relies on the following crucial fact \cite[Theorem 9.3]{BK2}; we state the fact for finite fields, and its analog for non-Archimedean local fields.

\begin{description}
 \item[\normalfont Over $\FF$ finite] The Braverman--Kazhdan gamma factor of an irreducible representation of a Deligne--Luztig representation attached to a character $\theta$ of  a Cartan subgroup $T_w$ of $\GL_n$ depends  only on the data $(T_w,\theta)$.
 \item[\normalfont Over $F$ local, non-Arch.] The Godement--Jacquet gamma factor of an irreducible subquotient of a representation of $\GL_n$ parabolically induced from a supercuspidal representation $\sigma$ of a Levi $L$ depends  only on the data $(L,\sigma)$.
 \end{description}

Let us now describe the construction  of the sheaf $\mathcal J_r$.  It uses the invariant-theoretic quotient $G\xrightarrow{p} \Dfrac{G}{G} \simeq T\sslash W$ (where $G$ acts by conjugation  on itself, and $T$ is its abstract Cartan). If $\iota: G^\rs\hookrightarrow G$ is the embedding of the regular semisimple part, then the sheaf $\mathcal J_r$ has the form
\[ \mathcal J_r = \iota_{!*} (p\iota)^* \mathcal J_r^{T\sslash W},\]
where $\mathcal J_r^{T\sslash W}$ is a certain sheaf on $T^\rs/W$, obtained from the analogous gamma-sheaf $\mathcal J_r^T$ for $T$. 

The gamma-sheaf $\mathcal J_r^T$ for $T$ is easy to describe; we keep assuming the existence of a central $\Gm\to \check G$ that composed with $r:\check G \to \GL(V)$ gives the scaling on $V$. The restriction $r|_{\check T}$ gives rise, dually, to a map $\pi: T_V\to T$, where $T_V$ is the Cartan of the Langlands dual (over $\FF$) of $\GL(V)$. There is a trace map $T_V \to \Ga$. We fix an Artin--Schreier sheaf $\mathcal L_\psi$ on $\Ga$, by choosing a nontrivial additive character $\FF\to \QQ_l^\times$, and take $\mathcal J_r^T = \pi_! \tr^* \mathcal L_\psi$. A natural $W$-equivariant structure, described in \cite[6.3]{Ngo-Takagi}, allows us to descend this sheaf from the regular set $T^\rs$, where $W$ acts freely, to the quotient $T^\rs/W$, obtaining the sheaf $\mathcal J_r^{T\sslash W}$. Finally, this sheaf, and hence the sheaf $\mathcal J_r$ on $G$, carry a natural Weil structure; thus, we can take their Frobenius trace.

\begin{theorem}[Braverman--Kazhdan, Cheng--Ngo, Chen]
 Let $\gimel_r$ denote the Frobenius trace of the sheaf $\mathcal J_r$ constructed above, multiplied by counting measure on $G(\mathbb F)$. It is a conjugation-invariant measure on $G(\mathbb F)$, and acts on every irreducible representation $\pi$ of $G(\mathbb F)$ by the gamma factor $\gamma(\pi,r)$.
\end{theorem}

A different proof of this result has appeared in a preprint of Laumon and Letellier \cite{Laumon-Letellier}. In a more recent preprint \cite{Laumon-Letellier2}, the same authors focus on the case of $r=$ symmetric-square representation of $G=\SL_2$, $\PGL_2$ or $\GL_2$, and investigate the question of whether the $r$-Fourier transform (composed with the Chevalley involution) extends to an involutive automorphism of the space of functions on the $\FF$-points of a stacky embedding of $G$.

\subsection{The case of local fields}

Following the resolution of this conjecture for finite fields, it is natural to ask for a description of the $r$-Fourier kernel for local fields. Some speculative proposals for this $r$-Fourier kernel were proposed in \cite[\S~7]{BK2} and  \cite[\S~6.2]{Ngo-Takagi}, modelled after the finite-field paradigm of obtaining the $r$-Fourier kernel on $G$ from the analogous $r$-Fourier kernel on $T$. These proposals should be viewed as motivational; as we will see, some correction is due, to make the kernel compatible with parabolic induction.
This correction was calculated by L.\ Lafforgue for the case of $\GL_2$. There is ongoing work by Ng\^o and Z.\ Luo, generalizing this to $\GL_n$, for arbitrary $n$. 

Let us discuss the solution proposed by Lafforgue, which appears in \cite[Expos\'e III]{Lafforgue-esquisse}. I hasten to clarify that some conventions will be different here than in \cite{Lafforgue-esquisse}; namely, all our kernels will be represented as \emph{a function multiplied by a Haar measure}, and I will use the same notation for the function and the distribution; that, of course, introduces a scalar ambiguity, which I will ignore, in this expository article. The reader can consult the original reference for details, but should be careful to notice that the measures used there, both on the reductive group and on its tori, are eigenmeasures that are not invariant.

Let $r=r_k=$ be the $k$-th symmetric power of the standard representation of $\GL_2$. This is actually a faithful representation of the group $\check G = \GL_2/\mu_k$, and we will construct the $r$-Fourier kernel on the dual group $G=\GL_2\times_{\Gm} \Gm$, where the map $\GL_2\to\Gm$ is the determinant, and the map $\Gm\to\Gm$ is the $k$-th power -- i.e., $G$ classifies automorphisms of a $2$-dimensional vector space, together with a $k$-th root of the determinant.

Recall the canonical identification $\Dfrac{G}{G}\simeq T\sslash W$. 
We will start by defining the ``naive'' kernel $\gimel_r^{T\sslash W}$ as a function on $(T\sslash W)(F)$ (or rather, its regular semisimple points), following \cite[\S~II.4]{Lafforgue-esquisse} and \cite[\S~6.2]{Ngo-Takagi}. To be clear, only its restriction to the image of $T'(F)\to (T\sslash W)(F)$, for $T'$ a split maximal torus, matters for the statement of ``compatibility with parabolic induction'' below, but presenting the entire kernel provides food for thought. 

The definition of $\gimel_r^{T\sslash W}$ is that it is that function on $T^\rs/ W$ (that is, on its $F$-points) whose pullback to the regular points of any Cartan subgroup $T'\subset G$ is the $r$-Fourier kernel $\gimel_r^{T'}$. Let us first explain this: The ``baby case'' is the Fourier kernel $|x|^\frac{1}{2} \psi(x) d^\times x$ for the action of $F^\times$ on $\mathcal D(F)$. Now, for any maximal torus $T'\subset G$, we have an embedding, up to conjugation, of its dual torus $\check T'\hookrightarrow \check G$. Restricting $r$ to $\check T'$ and diagonalizing, we obtain a map $\check T'\to \Gm^{k+1}$, whose dual is a map ${\Gm}_{\bar F}^{k+1} \to T'_{\bar F}$ (i.e., defined over the algebraic closure). The Galois group $\Gamma$ acts on $\check T'$ and on the set of weights of $r|_{\check T'}$ (which are all of multiplicity one in this case -- this simplifies things a little), hence this map descends to a map of tori $p:T_E'\to T'$, where $T_E' = \Res_{E/F}\Gm$ is the restriction of scalars of $\Gm$ from some separable $F$-algebra $E$. On the other hand, $T_E'$ comes with a ``trace'' map to $\Ga$ (descending from the sum of coordinates on ${\Gm}_{\bar F}^{k+1}$). We then take 
\[\gimel_r^{T'} = p_* \left(|\bullet|_E^\frac{1}{2} \tr^* \psi\right),\] 
i.e., we pull back the additive character $\psi$ from $F$, multiply it by the square root of the modulus character for the action of $T_E'(F)\simeq E^\times$ on $E$, and push it forward (as a measure, after multiplying by a Haar measure) to $T'(F)$. As mentioned, we will not bother to keep track of Haar measures here, but there is a standard way to choose them compatibly for all tori. 

We choose coordinates $(c,a)\in \Ga\times\Gm$ for the quotient $\Dfrac{G}{G}= T\sslash W$, where we send $(g,a)\in \GL_2\times_{\Gm} \Gm$ to $c=\tr g$ and $a$.
Having defined the naive kernel $\gimel_r^{T\sslash W}$, Lafforgue corrects it to a kernel $\gimel_r^G$ given by the formula
\[ \mathcal F_c \gimel_r^G(\xi,a) = |\xi|\cdot |a|^\frac{k}{2} \mathcal F_c  \gimel_r^{T\sslash W}(\xi,a),\]
where $\mathcal F_c$ is the usual Fourier transform in the $c$-variable (with $\xi$ denoting the dual variable).
I do not know a conceptual explanation for this formula, but the correction is necessary for compatibility with parabolic induction. Namely, in \cite[Proposition III.2]{Lafforgue-esquisse}, Lafforgue performs a formal calculation that suggests the following property (which fails for the naive kernel):
\[ (\gimel_r^G \star \Phi)_N = \gimel_r^T \star \Phi_N,\]
for all Schwartz functions $\Phi$ on $G$, where $\Phi_N$ is the normalized constant term
\[ \Phi_N(t) = |\delta_B(t)|^{-\frac{1}{2}} \int_N \Phi(nt) dn.\]
More conceptually stated, the action of $\gimel_r^G$ on half-densities on $N\backslash G$ matches the action of $\gimel_r^T$ under the ``left'' action of $T$ on $N\backslash G$. This is necessary, for $\gimel_r^G$ to act by the correct gamma factors on principal series representations.

\subsection{Fourier transforms on other nonlinear spaces} \label{ss:Fourier-nonlinear}

There are more conjectural nonlinear Fourier transforms than the $r$-Fourier transforms in the case of the group. In fact, the local functional equations of various classical integral representations of $L$-functions can be construed as some kind of ``Fourier transforms.''

Take, for example, Hecke's classical integral representation for the standard $L$-function of a modular form, as reinterpreted in the adelic language by Jacquet and Langlands \cite{JL}. Hence, let $V_2=V_1\oplus V_1'$ be a $2$-dimensional vector space written as a direct sum of two $1$-dimensional subspaces, and let $X=\GL(V_1)\backslash \GL(V_2)$, $X^\vee =\GL(V_1')\backslash \GL(V_2)$. If $w\in G=\GL(V_2)$ is an element interchanging the subspaces $V_1$ and $V_1'$, it defines a $G$-equivariant isomorphism $X\xrightarrow\sim X^\vee$, sending the coset of $g$ to the coset of $wg$. This gives rise to an isometric isomorphism of Schwartz spaces of functions or half-densities,
\[ \mathcal D(X) \xrightarrow\sim \mathcal D(X^\vee),\]
which is the basis of the local functional equation of the standard $L$-function, see \cite{JPSS}. To formulate the exact functional equation with the right gamma factors, one needs to explain how to choose Whittaker models compatibly with the choice of $w$, a topic that deserves close attention, but not in the current exposition. I also point the reader to the recent thesis of G.\ Dor \cite{Dor}, for a very novel and deep interpretation of the relationship between the Godement--Jacquet and Jacquet--Langlands (Hecke) integral representations (without, however, a comparison of local gamma factors). 

For an example that looks more like a classical Fourier transform, consider the Rankin--Selberg integral for $G=\GL_n\times_{\det}\GL_n$, which can be written as the integral of an automorphic form on $G$ against a theta series coming from the space of Schwartz functions on the ``Rankin--Selberg variety'' $X=\Std\times^{\GL_n^\diag} G$, where $\Std$ denotes the vector space of the standard representation. (See the article of Beuzart-Plessis, \cite{Beuzart-Plessis} for the definition of theta series.) A ``dual'' space is obtained by considering the dual representation, $X^\vee = \Std^\vee \times^{\GL_n^\diag} G$, and fixing an additive character we get an equivariant Fourier transform of Schwartz half-densities (defined along the dual fibers of the maps $X\to \GL_n^\diag\backslash G \leftarrow X^\vee$, that is, their points over a local field):
\[ \mathcal D(X)\xrightarrow\sim \mathcal \mathcal D(X^\vee).\]
This Fourier transform gives rise to the functional equation of the Rankin--Selberg $L$-function.

More generally, let $X$ be any affine spherical $G$-variety, for a quasisplit group $G$, possibly with some additional restrictions (e.g., without ``spherical roots of type $N$''). Let $X^\vee$ denote the same $G$-variety, but with the $G$-action composed with a Chevalley involution. ``A Chevalley involution'' is an involution of $G$ that acts as minus the longest element of the Weyl group on the (abstract) root system of $G$. All such involutions are conjugate over the algebraic closure, but the choice of its $G(F)$-conjugacy class may matter for what follows. For some groups (but not all, e.g., not for $\SL_n$), there is a ``duality involution'' that satisfies $\pi^c \simeq \tilde\pi$; I point the reader to the relevant article of D.\ Prasad \cite{Prasad-contragredient}.  In any case, a Chevalley involution $c$ is expected to interchange the $L$-packets of an irreducible representation $\pi$ and its contragredient $\tilde\pi$. 

In the setting above, one often expects to be able to define a ``Fourier transform,'' which is again an $L^2$-isometric, $G$-equivariant isomorphism between spaces of Schwartz half-densities on $X$ and $X^\vee$.  As we saw in Section \ref{sec:Schwartz}, when the affine variety $X$ is singular, the appropriate definition of Schwartz spaces is still unclear. These uncertainties notwithstanding, we should highlight a number of promising cases in the study of such nonlinear Fourier transforms:

\subsubsection{The parabolic case} This case is probably the origin of all ideas discussed in this section. At the level of functions, it is discussed again in a couple of other papers of Braverman and Kazhdan \cite{BK1, BK3}, but it is closely related to ideas of Laumon and Drinfeld on ``compactified'' Eisenstein functors in the geometric setting \cite{BG}. In order to stay closer to the existing literature, and to make connections to the theory of Eisenstein series, let us work with functions instead of densities here, using the letter $\mathcal S$ to denote Schwartz spaces of functions.  

The idea here is to attach a suitable Schwartz space of functions $\mathcal S(\overline{U\backslash G})\subset C^\infty(U\backslash G)$, where $U\subset P\subset G$ is the unipotent radical of a parabolic subgroup. In principle, this space is not associated to the homogeneous space $U\backslash G$, but to its affine closure $\overline{U\backslash G} = \Spec F[U\backslash G]$. 

The ``baby case'' of this is as classical as Eisenstein series themselves: It corresponds to the difference between the definitions
\[ E(z,s) = \sum_{(m,n)=1} \frac{y^s}{|mz+n|^{2s}}\]
and
\[ E^*(z,s) = \pi^{-s} \Gamma(s) \sum_{(m,n)\ne(0,0)} \frac{y^s}{|mz+n|^{2s}}.\]
The two are related by $E^*(z,s) = Z(2s) E(z,s)$ (where $Z(s)$ denotes the complete zeta function),
making the functional equation for $E^*(z,s)$ nicer:
\[ E^*(z,s) = E^*(1-z.s).\]

Adelically, the ``normalized'' Eisenstein series $E^*(z,s)$ can be obtained from the Schwartz space of the affine plane $\mathbb A^2 = \overline{U\backslash\SL_2}$ by 
\[ \mathcal S(\mathbb A_k^2) \to I_B^G(\chi) \xrightarrow{\mathcal E_\chi} C^\infty([G]),\]
where $\mathbb A_k$ denotes the adeles of a global field $k$,   $I_B^G(\chi)$ denotes the principal series representation unitarily induced from a suitable character $\chi$ (corresponding to the parameter $s$ above), and $\mathcal E_\chi$ is the Eisenstein series morphism to functions on $[G]=G(k)\backslash G(\mathbb A_k)$. This is to be juxtapposed to the ``usual'' Eisenstein series $E(z,s)$, which comes from the similar diagram by replacing $\mathcal S(\mathbb A_k^2)$ by $\mathcal S(U\backslash G(\mathbb A_k))$. The functional equation for $E^*(z,s)$ is a corollary of the Poisson summation formula for the lattice $k^2\subset \mathbb A_k^2$.

In general, however, the affine closure $\overline{U\backslash G}$ is a singular variety, and does not have a vector space structure. Nonetheless, it was observed in the geometric Langlands program that Eisenstein series defined by ``Schwartz sheaves'' on such spaces are better behaved. Braverman and Kazhdan \cite{BK1, BK3} applied this wisdom to the function-theoretic setting, where they managed to provide a definition for the Schwartz spaces $\mathcal S(\overline{U\backslash G})$ when $U$ is the unipotent radical of a Borel subgroup, and more generally $\mathcal S(\overline{[P,P]\backslash G})$, for any parabolic $P$. They also defined Fourier transforms (normalized intertwining operators) of the form
\[ \mathcal S(\overline{[P,P]\backslash G}) \xrightarrow\sim \mathcal S(\overline{[P^-,P^-]\backslash G}),
\]
where $P^-$ is an opposite parabolic to $P$, and proved a Poisson summation formula (eventually, by reducing it to the case of $\SL_2$ discussed above), giving rise to the functional equation of the normalized Eisenstein series. (A more complete version of this Poisson summation formula was proven by Getz--Liu \cite{GL-Poisson}.)

For the general case $X=\overline{U\backslash G}$ (that is, when $P$ is not a Borel subgroup), it is worth noting that \emph{even in the geometric Langlands program, the conjectural Fourier transform is absent}; undoubtedly a deep problem, whose resolution would be very interesting. We do know, however, from the calculations of intersection complexes in  \cite{BFGM}, that the IC function of $X$ gives rise to \emph{normalized Eisenstein series}; more precisely, the function-theoretic interpretation of \cite[Theorem 1.12]{BFGM} is 
\[ IC_{\overline{U\backslash G}} = L(\check{\mathfrak u},1) \star 1_{U\backslash G(\mathfrak o)},\]
where the notation is as follows:
\begin{itemize}
 \item[--] $\mathfrak o = \FF_q[[t]]$ is a DVR in equal characteristic, whose fraction field we will denote below by $F$;
 \item[--] $1_{U\backslash G(\mathfrak o)}$, the basic function of $U\backslash G$, is the characteristic function of its integers;
 \item[--] $IC_{\overline{U\backslash G}}$ is the basic function  of $\overline{U\backslash G}$, define via the intersection complex of its arc space as in \cite{BNS}; it is a smooth function on $U\backslash G(F)$, whose support belongs to $\overline{U\backslash G}(\mathfrak o)$;
 \item[--] finally, $L(\check{\mathfrak u},1)$ is a series of elements of the unramified Hecke algebra of the Levi quotient $L$ of $P$; it is the series that, under the Satake isomorphism, corresponds to $\sum_{i\ge 0} q^{-i} \Sym^i \check{\mathfrak u}$, where $\check{\mathfrak u}$ is the representation of the dual Levi $\check L$ on the dual nilpotent radical $\check{\mathfrak u}$. The convolution action (denoted by $\star$) of the Hecke algebra of $L$ on functions on $U\backslash G$ is obtained from the $L^2$-normalized action of $L$ (i.e., we are secretly working with half-densities); the reader can pin down the details by keeping in mind that, here, we define the action of $L$ on $U\backslash G$ as a left action, and for $G=\SL_2$, $L(\check{\mathfrak u},1) \star 1_{U\backslash G(\mathfrak o)}$ stands for the characteristic function of $\mathfrak o^2$. 
\end{itemize}

This suggests that the normalized Fourier transforms 
\[\mathcal S(\overline{U\backslash G}) \xrightarrow\sim \mathcal S(\overline{U^-\backslash G}) \]
should spectrally decompose in terms of the \emph{normalized intertwining operators} of Langlands and Shahidi \cite{Langlands-Eisenstein, Shahidi-Plancherel}. In the case of degenerate Eisenstein series (that is, replacing $U$ by $[P,P]$), for classical groups, this was verified by Shahidi (and W.\ W.\ Li) in \cite{Shahidi-generalizedFourier}. A direct proof of the Poisson summation formula for these Fourier transforms, leading to the meromorphic continuation of the normalized Eisenstein series, would greatly simplify and generalize the Langlands--Shahidi method, as it would imply the meromorphic continuation of the ratio between normalized and unnormalized Eisenstein series.

\subsubsection{Other cases}

There are several other affine spherical varieties $X$ which turn out to be fiber bundles of the form $X\to H\backslash G$, with $H$ reductive and fiber isomorphic to $Y'=\overline{[P',P']\backslash H'}$, for some parabolic $P'$ in a reductive group $H'$ containing $H$. This explains a large part of the Rankin--Selberg method, as was observed in \cite{SaRS}, with the functional equation obtained from Fourier transform on $Y'$. For example, the affine closure of the quotient of $\SL_2^2$ by the subgroup of triples of upper triangular matrices with entries satisfying $x_1+x_2+x_3=0$ happens to be isomorphic to the affine closure of $[S,S]\backslash\Sp_6$, where $S$ denotes the Siegel parabolic, and this fact is behind Garrett's integral representation of the triple product $L$-function \cite{Garrett}. 

An entirely new case of nonlinear Fourier transforms, on a space which does not belong to the family of spaces such as $Y'$ above, was considered by Getz and Liu in \cite{GL-triples} (and refined in \cite{GH, GHL}. In it, they consider a triple $(V_i, q_i)_{i=1}^3$ of even-dimensional (nondegenerate) quadratic spaces, and consider the affine variety $X$ given by the equation 
\[ q_1 = q_2 = q_3.\]
This variety is spherical under the action of $G=$ the fiber product of orthogonal similitude groups $\GO(V_1)\times_{\Gm}\GO(V_2)\times_{\Gm}\GO(V_3)$. In the series of aforementioned papers, for $F$ a local field in characteristic zero, a space $\mathcal S(X)$ of smooth functions on $X^\infty(F)$ is defined (where $X^\infty\subset X$ is the smooth locus), as the image of a morphism 
\[ \mathcal S(Y)\otimes \mathcal S(V) \to \mathcal S(X),\]
where the notation is as follows: $V = \bigoplus_i V_i$, and $Y=$ the affine closure of $[S,S]\backslash \Sp_6$, with notation as above. The construction of the morphism is via the Weil representation for $\SL_2^3 \times G$, with $\SL_2^3$ considered as a subgroup of $\Sp_6$, and the morphism is $\SL_2^3$-invariant, hence factors through the coinvariant quotient 
\[ \left(\mathcal S(Y)\otimes \mathcal S(V) \right)_{\SL_2^3} \to \mathcal S(X).\]
Although it is not proven in these papers, it may help, to fix ideas, to think of $\mathcal S(X)$ as an avatar for this coinvariant space.

Then, the authors of \cite{GH} define a Fourier transform on $\mathcal S(X)$ by descending the Braverman--Kazhdan Fourier transform of $\mathcal S(Y)$. Note that $Y^\vee\simeq Y$ under the action of $\Sp_6$, although this inverts the commuting action of $\Gm=S/[S,S]$; for uniformity with the setup above, we should be thinking of the Fourier transform of Getz et.~al.\ as a $G$-equivariant isomorphism $\mathcal F_X: \mathcal S(X)\xrightarrow\sim \mathcal S(X^\vee)$. In \cite{GHL} it is shown that this transform extends to an $L^2$-isometry. Finally, when all spaces are defined over a number field $k$, a Poisson summation formula is proven in \cite{GL-triples, GH}, which has the form
\[ \sum_{\gamma\in Y^\bullet(k)} \Phi(\gamma) + \mbox{(boundary terms)} 
 = \sum_{\gamma\in Y^\bullet(k)} \mathcal F_X \Phi(\gamma) + \mbox{(boundary terms)}, 
\]
where $Y^\bullet$ is the open $G$-orbit. The expression of the boundary terms is not intrinsic to $\Phi$ at this point, but at least one can impose local restrictions on $\Phi$ that will cause them to vanish.

\subsection{Hankel transforms for the (relative) trace formula} \label{ss:Hankel}

\subsubsection{Hankel transforms for the standard relative trace formula}
The Fourier transforms $\mathcal D(X) \xrightarrow\sim \mathcal D(X^\vee)$ that we speculated on in the previous subsection are likely not available for every space $X$, and there are reasons for that on the side of $L$-functions: It is sometimes the square of the absolute value of a period that is related to an $L$-value, not the period itself, and there is no reason to expect the ``square root'' represented by the period to carry some sort of ``functional equation.'' In those cases, one should have the appropriate Fourier transform appearing only at the level of the quotients of the relative trace formula, which spectrally encodes this square of the period.

Rather than speculating in general, let me explain this idea in a setting that has already been worked out by H.\ Xue in \cite{Xue} (with the lowest-rank case appearing previously in \cite{SaBE2}): Let $E/F$ be a quadratic extension of local fields (in characteristic zero, in the references that we will cite), and $A$ a central simple algebra over $F$ of dimension $4n^2$, together with an embedding $E\to A$. Let $B$ be the centralizer of $E$ in $A$. Let $G=A^\times/\Gm$, $H=B^\times/\Gm$, considered as algebraic groups over $F$. When $A=\Mat_2(F)$, $H$ is simply a maximal torus in $G=\PGL_2$, and we are in the setting of Waldspurger's periods \cite{Waldspurger}. The dual group of the variety $X=H\backslash G$ is $\Sp_{2n}$, and it is a general expectation  that, in the same setting over global fields, for non-residual automorphic representations $\pi$ on $G$ whose exterior-square $L$-function has a pole at $s=1$ (indicating a global lift from $\Sp_{2n}$), and 
satisfying appropriate local conditions, one has an Ichino--Ikeda-type formula decomposing the square of the absolute value of the $H$-period, 
\[ \pi \ni \phi \mapsto \left|\int_{[H]} \phi(h) dh\right|^2\]
into an Euler product of local functionals. (I do not know of such a precise factorization in the literature, but see \cite{GJ, FMW} for weaker versions of the conjecture and global results; one could, more generally, integrate against an automorphic character of $H$, but we will restrict ourselves to the trivial character, both globally and locally.) These local functionals should evaluate, at almost every place $F$, to the quotient of local $L$-factors
\[\frac{L(\pi_E, \frac{1}{2}) L(\pi_F, \mathfrak{sl}_{2n}/\mathfrak{sp}_{2n},1)}{L(\pi_F, \mathfrak{sp}_{2n},1)},\] 
where $\pi_E$ is the base change of $\pi_F$ to $E$. 

The local counterpart of this conjecture, providing the local conditions for it, is a conjecture of Prasad--Takloo-Bighash \cite{PTB}, generalizing a theorem of Tunnell and Saito \cite{Tunnell, Saito}, which states that a tempered irreducible representation $\pi$ of $G$ is $H$-distinguished if and only if its Langlands parameter factors through $\Sp_{2n} \to \GL_{2n}$, and $\epsilon(\pi_E)\eta(-1)^n=(-1)^r$, where $\epsilon(\pi_E)$ is the local root number of the base change of $\pi$ to $E$, $\eta$ is the quadratic character attached to $E/F$, and $(r-1)$ is the split rank of $G$. This conjecture has been proven for discrete series by H.\ Xue and M.\ Suzuki \cite{Xue, SX}; at the heart of the argument in \cite{Xue} is an involution on the space $\mathcal S(H\backslash G/H)$ of test measures for the relative trace formula of the quotient $H\backslash G/H = (X\times X)/G$. We will think of this involution as the analog of the Fourier transform on the space of functions or half-densities on $X$, and will call it a ``Hankel transform.'' (That name has been used by Ng\^o to describe the $r$-Fourier transforms discussed previously, but we will reserve this name for analogous transforms at the level of trace formulas.)

The space of test measures $\mathcal S(H\backslash G/H)$ should be understood of as the direct sum, over isomorphism classes of central simple algebras of dimension $4n^2$ with an embedding $E\hookrightarrow A$ as above, of the corresponding coinvariant spaces $\mathcal S(H\backslash G)_H$. The said involution acts by $(-1)^r$ on the summands with split rank $(r-1)$, and on the other hand Xue shows that it acts by a factor of $\epsilon(\pi_E)\eta(-1)^n$ on the relative character of $\pi$. It is quite notable that there is a meaningful, nontrivial version of Fourier transform in this setting, where $X^\vee \simeq X$ and the associated $L$-function only appears by its value at the central point $\frac{1}{2}$; still, its functional equation is expressed in a nontrivial way by the epsilon factor. (It is best to think of the ``Hankel transform'' here as being the product of Xue's involution by $\eta(-1)^n$.)

\subsubsection{Hankel transforms for the (relative) trace formula with nonstandard test measures} \label{sssHankel}

Let us place ourselves in the conjectural setting of the Braverman--Kazhdan--Ng\^o program, with a nonstandard Schwartz space $\mathcal S_r(G)$ associated to a representation $r$ of the $L$-group of $G$, together with an $r$-Fourier transform (which, for now, we take to be between spaces of measures)
\[ \mathcal F_r: \mathcal S_r(G) \xrightarrow\sim \mathcal S_{r^\vee}(G).\]

Descending to test measures for the Arthur--Selberg or the Kuznetsov trace formulas, this would give rise to ``Hankel transforms''
\[ \mathcal H_r: \mathcal S_r(\frac{G}{G}) \xrightarrow\sim \mathcal S_{r^\vee}(\frac{G}{G}),\]
or
\[ \mathcal H_r: \mathcal S_r((N,\psi)\backslash G/(N,\psi)) \xrightarrow\sim \mathcal S_{r^\vee}((N,\psi)\backslash G/(N,\psi)),\]
between nonstandard spaces of test measures for the corresponding trace formula. If we could prove a Poisson summation formula for these transforms, in the sense that, for the corresponding spaces over the adeles, the diagram 
\[\xymatrix{ \mathcal S_r((N,\psi)\backslash G/(N,\psi)) \ar[rr]^{\mathcal H_r}\ar[dr]_{\KTF} & & \ar[dl]^{\KTF}\mathcal S_{r^\vee}((N,\psi)\backslash G/(N,\psi))\\
& \CC &}\]
commutes, that would give rise to a proof of the functional equation of the $L$-function corresponding to $r$. To be clear, there are serious difficulties making sense of the ``trace formula'' functional on nonstandard spaces of test measures, which we will discuss in Section \ref{sec:global} below.

Such transforms have been computed, in the setting of the Kuznetsov formula in a few cases, namely:
\begin{enumerate}
 \item For $r=$ the standard representation of $\GL_n$, by Jacquet \cite{Jacquet-smoothtransfer}. In this case, the nonstandard space of test functions and the Hankel transform descend from the Schwartz space of $\Mat_n$, and its Fourier transform.
 \item For $r=$ the symmetric-square representation when $G=\SL_2$, by the author \cite{SaTransfer2}.
\end{enumerate}

A formula for the case of  $r= \Std \oplus \Std \otimes \eta$, $G=\PGL_2$, with $\eta$ a (possibly trivial) quadratic character of the Galois group,  was also proven in \cite{SaBE2}. However, this case can easily be reduced to repeated applications of the first case above.

A remarkable feature of these formulas is that, despite the non-abelian nature of the $L$-functions whose functional equation they encode, they are given by a combination of abelian Fourier transforms and various ``correction factors'' -- much like the transfer operators that we encountered in Sections \ref{sec:localtori}, \ref{sec:samerank}. For example, let us discuss the case of the standard representation of $\GL_2$. It turns out that the formula is nicer if we work with half-densities, instead of measures (see \cite[Theorem 8.1]{SaHanoi} for details on how to translate Jacquet's result to half-densities). Hence, varying our discussion in \S~\ref{ss:Kuztori}, we will introduce a nonstandard space $\mathcal D_{\Std}=\mathcal D_{\Std}((N,\psi)\backslash G/(N,\psi))$ of ``test half-densities'' for the Kuznetsov formula, considered as smooth half-densities on the torus $T$ of diagonal elements via the embedding $t\mapsto w t$, with $w$ the permutation matrix, and defined by 
the formula 
\[ f(t)= (|\zeta| d\zeta)^\frac{1}{2} \cdot \int_{F^2} \Phi \left(\begin{pmatrix} 1& x \\ &1\end{pmatrix} wt\begin{pmatrix} 1& y \\&1 \end{pmatrix} \right) \psi^{-1}(x+y) dx dy,\]
where $\Phi$ ranges over Schwartz functions on $\Mat_2$. We also define the analogous space $\mathcal D_{\Std^\vee}$, by using the natural embedding of $\GL_2$ into the dual vector space. There is then a Hankel transform $\mathcal H_{\Std}: \mathcal D_{\Std} \to \mathcal D_{\Std^\vee}$, descending from the equivariant Fourier transform $\mathcal D(\Mat_2)\to \mathcal D(\Mat_2^\vee)$, and Jacquet's formula reads
 \begin{equation}\label{Hankel-Jacquet}
  \mathcal H_{\Std} f = D_{-\check\epsilon_1,\frac{1}{2}} \star \left( \psi^{-1}(e^{-\alpha}) \circ D_{-\check\epsilon_2,\frac{1}{2}} \star f\right),
 \end{equation}
where our notation is as in \S~\ref{ss:KuzSLn}. Namely, $D_{-\check\epsilon_i,s}$ is convolution on $T$ with the measure $D_s$ (see \eqref{Ds}), pushed forward to $T$ by the cocharacter $(-\check\epsilon_i)$, where $\epsilon_i$ is the standard cocharacter into the $i$-th coordinate of $T$.

The aforementioned ``correction  factor'' here is  $\psi^{-1}(e^{-\alpha})$, which simply means multiplication by the function $\begin{pmatrix} a &  \\ &d \end{pmatrix}\mapsto \psi^{-1}(\frac{d}{a})$. 
These correction factors are poorly understood; it would be desirable to have an interpretation of these transfer operators in terms of quantization, as for the transfer operators discussed in \S~\ref{ss:transferrankone}.
I will refrain from presenting these formulas here, pointing the reader to the expository article \cite{SaHanoi}.

\section{Global comparisons}\label{sec:global} 

Beyond Endoscopy is a strategy to solve the global problem of functoriality for automorphic representations, and some of the most involved work in this context has been in the global setting. It is quite a challenge, however, to summarize and put in a common context the global results that have appeared so far in the literature, and the techniques developed to achieve them. In order to do so, I will focus on aspects that can be understood in combination with the local results that were mentioned so far, although the relation between global and local aspects has not been sufficiently explored in the literature, and will be slighty speculative. As a result, the present article is purporting to be neither comprehensive, nor loyal to the point of view of the original authors. 

We will focus on 4 pieces of work, in reverse chronological order: The work of Frenkel--Langlands--Ng\^o and A.\ Altu\u{g} on the trace formula for $\GL_2$, the author's work on the relative trace formula for torus periods on $\PGL_2$, the work of E.\ Herman on the functional equation of the standard $L$-function for $\GL_2$-automorphic forms, and the work of A.\ Venkatesh on stable functoriality from tori to $\GL_2$.

\subsection{Poisson summation for the trace formula of $\GL_2$} \label{ss:Altug}

One of the main goals of the work of Frenkel--Langlands--Ng\^o and A.\ Altu\u{g} \cite{FLN, Altug1, Altug2, Altug3} on the trace formula for $\GL_2$ was to extract from it an expression that contains only the ``Ramanujan'' part of the spectrum, removing the contributions of residual automorphic representations (characters, in the case of $\GL_2$). 

To simplify notation, we will work with the stable trace formula of $\SL_2$, so that the only $1$-dimensional automorphic representation is the trivial one, and the Steinberg--Hitchin base $\mathfrak C=\Dfrac{G}{G}$ is isomorphic to $\mathbb A^1$ via the trace function. Every mention of ``orbital integrals'' and ``trace formula''  below will refer to stable orbital integrals, and the stable trace formula.

The fundamental observation  of \cite{FLN} is that the trivial representation should correspond to the value at $0$ of the ``Fourier transform of orbital integrals.'' This is, first of all, a local statement, and it is clear if we think of test measures instead of test functions and their orbital integrals: Considering test spaces of measures, the composition of pushforwards under $G\to \Dfrac{G}{G} \to *$ is the integral of a measure -- i.e., the trace of its action on the trivial representation -- and at the same time it is equal to the value at $0$ of the Fourier transform of its pushforward to $\mathfrak C=\Dfrac{G}{G}$.

In order to express the global trace formula, however, in terms of the Fourier transform, one needs a version of Poisson summation formula. Namely, very simplistically thinking, we can think of the geometric side of the trace formula as a sum over the rational points of $\mathfrak C$, and we need to convert it to a sum over ``the dual vector space'' $\mathfrak C^*$. (Note, also, the conceptual mystery of the fact that we treat $\mathfrak C$ as a vector space; there is not a clear reason why we should do so.) In reality, of course, this sum only applies to the elliptic terms in the trace formula, with the actual expression being more complicated, depending on the version of the trace formula that  one chooses to use (e.g., non-invariant or invariant).  This is only part of the difficulty; another fundamental difficulty is that orbital integrals/test measures do not belong to the standard Schwartz space of functions/measures on $\mathbb A^1$, hence are not immediately suited for an application of the Poisson summation formula. 

What form should one expect the trace formula to take, after application of the Poisson summation formula? Altu\u{g}'s work made the biggest progress to answer this question, but before we get to it, we will reformulate the question, in view of the local theory already developed. Namely, we saw in \S~\ref{ss:KuzSL2} that Fourier transform on $\mathfrak C$, locally, has a special meaning: It is a transfer operator between stable test measures for the trace formula of $\SL_2$, and test measures for its Kuznetsov formula. From this point of view, the sought Poisson summation formula would, roughtly, correspond to the statement that the following diagram commutes:
\begin{equation}\label{STFKTF}\xymatrix{ \mathcal S(\frac{G}{G})^\st \ar[rr]^{\mathcal T^{-1}}\ar[dr]_{\STF} & & \ar[dl]^{\KTF}\mathcal S_{L(\Ad,1)}((N,\psi)\backslash G/(N,\psi))\\
& \mathbb C &}\end{equation}
where the notation is as follows:
\begin{itemize}
 \item[--]  $\mathcal S(\frac{G}{G})^\st $ is the stable space of test measures for the trace formula, over the adeles, defined as a restricted tensor product of the local spaces;
 \item[--] $\mathcal T^{-1}$ is the Fourier transform, which corresponds to the inverse of the transfer operator discussed in \S~\ref{ss:KuzSL2}, and $\mathcal S_{L(\Ad,1)}((N,\psi)\backslash G/(N,\psi))$ is its image, a nonstandard space of test measures for the Kuznetsov formula. This is also a restricted tensor product over all places, but with respect to a nonstandard ``basic vector,'' equal to the image of the basic vector of $\mathcal S(\frac{G}{G})^\st$ under $\mathcal T^{-1}$. (We will have more to say about this nonstandard basic vector below.)
 \item[--] ``STF'' stands for the functional of the stable trace formula, and ``KTF'' stands for the functional of the Kuznetsov formula.
\end{itemize}

The last point is not easy to state rigorously, since neither of the two functionals is well-defined. The stable trace formula that we have in mind is an idealistic ``naive invariant'' stable trace formula, where the contribution of a semisimple stable orbit represented by an element $\gamma$ would be the stable orbital integral at $\gamma$, multiplied by the volume of $G_\gamma(k)\backslash G_\gamma(\mathbb A_k)$ -- but this volume is infinite, for hyperbolic conjugacy classes. The Kuznetsov trace formula is well-defined on the standard space of test measures, but not on the nonstandard one, where the basic vector at almost every place $v$ is the generating series of the local unramified $L$-function $L(\pi_v, \Ad,1)$, as in \S~\ref{ss:KuzSL2}; the expression cannot be absolutely convergent, since on the spectral side the corresponding Euler products do not converge, and in fact $L(\pi, \Ad,s)$ has a pole at $s=1$ when $\pi$  is in the continuous (Eisenstein) spectrum. 

Nonetheless, we can view the work of Altu\u{g} as an attempt to make sense of the diagram \eqref{STFKTF}, or an approximation of it. This elucidates several aspects of his work, such as the appearance of Kloosterman sums in the expressions he obtains after Poisson summation -- such sums are inherent in the Kuznetsov formula. For fans of the Selberg trace formula this might be a disappointing point of view: already in his 2001 letter to Langlands \cite{Sarnak}, Sarnak pointed out that the Kuznetsov formula could be used to avoid the non-Ramanujan part of the spectrum, so one could say that we have not come a long way since then. Nonetheless, the global comparison of trace formulas, i.e., the effort to make sense of the Poisson summation formula of diagram \eqref{STFKTF}, has a lot to teach us about global techniques that will be necessary for the Beyond Endoscopy program. Moreover, as we stressed in \S~\ref{ss:KuzSL2}, such a comparison \emph{is} a special case of functoriality -- between different homogeneous spaces (the group $\SL_2$ and its Whittaker model), rather than different reductive groups.

Many of the interesting takeaways from Altu\u{g}'s work have been summarized by Arthur \cite{Arthur-BE}, who took them as a starting point to propose analogous questions for the trace formula in higher rank. I have little to add to the insights of the absolute expert in the field, and point the readers to his article. I will only attempt to highlight a couple of points which arise from thinking about Fourier transforms in terms of the comparison with the Kuznetsov formula:

Thinking about the spectral side of the Selberg trace formula and the Kuznetsov formula, two differences stand out: First, the traces of Hecke operators on cuspidal automorphic representations come weighted by the factor $\frac{1}{L(\pi,\Ad,1)}$ in the Kuznetsov formula, but this factor does not appear on the spectral side of the Selberg trace formula. Secondly, the Selberg trace formula (for $\SL_2$) has a contribution from the trivial representation, and also has a discrete contribution from Eisenstein series induced from characters of the Cartan that are fixed under the Weyl group; there is no analog of those terms on the spectral side of the Kuznetsov formula.

It is easy to explain how the first difference is resolved, in the putative comparison \eqref{STFKTF}: as already mentioned, the image of the basic vector, at almost every place $v$, under the transfer operator $\mathcal T^{-1}$ is the generating series of the local unramified $L$-factor $L(\pi_v, \Ad,1)$, hence the spectral contributions to the Kuznetsov formula with nonstandard test measures will come with an additional factor of $L(\pi, \Ad,1)$, canceling the difference between the 2 sides. (For Eisenstein series induced by a character $\chi$, it should be mentioned, the spectral contributions to the standard Kuznetsov formula come weighted by $\frac{1}{L(\chi,\check{\mathfrak g}/\check{\mathfrak t},1)}$, and multiplying them by $L(\pi, \Ad,1) = L(\chi,\check{\mathfrak g},1)$ introduces a pole; as we mentioned, this is the reason why the ``KTF'' functional is not well-defined for this space of nonstandard test measures, and similarly the contribution of Eisenstein series to the  ``naive invariant'' Selberg trace formula is infinite.)

For the second difference, we have explained the reasons that, locally, the contribution of the trivial representation is picked up by the value of the Fourier transform of a test measure at $0\in\mathfrak C^*$. As can be seen by the explicit form of the transfer operator in \eqref{inversetransfer}, the spaces $\mathfrak C^*\simeq \mathbb A^1$ and $N\backslash G\sslash N\simeq \mathbb A^1$ (the latter being the space of orbits for the standard Kuznetsov formula) cannot be identified; rather, there is a birational morphism $\xi\mapsto \xi^{-1}$ between them, so that the point $\xi=0\in \mathfrak C^*$ corresponds to ``infinity'' on the other space. This means that the contribution of the trivial representation is picked up by the behavior of the nonstandard test measures in $\mathcal S_{r}((N,\psi)\backslash G/(N,\psi))$ at ``infinity;'' note that the subspace 
$\mathcal S((N,\psi)\backslash G/(N,\psi))$ of standard test measures (locally) consists of measures that vanish there.

Globally, that means that, to achieve our putative comparison \eqref{STFKTF}, we need to include a contribution from the ``point at $\infty$'' of the space $N\backslash G\sslash N$, in a neighborhood of which our nonstandard test measures, unlike the standard ones, have nontrivial (nonzero) behavior. Vice versa, $\xi=\infty$ corresponds to the point of ``singular'' point of $N\backslash G\sslash N$ represented by the identity coset, and this is evidenced by the fact that the dual variable $\xi$ in the variant of the Kloosterman sums defined in \cite[Theorem 1.1]{Altug1} appears in the numerator, rather than the denominator, of the argument of the exponential -- compare with the variable $c$ in \cite[Proposition 3.7]{KLP}. (This by no means purports to be a full explanation of the idea that Altu\u{g}'s expression is related to the Kuznetsov formula; I am just highlighting some similarities, that the interested reader can investigate further.)

However, as is also highlighted by Arthur in \cite{Arthur-BE}, the point $\xi=0\in\mathfrak C^*$ in Altu\u{g}'s work contains not only the contribution of the trivial representation, but also the discrete contribution of the Eisenstein series unitarily induced from quadratic characters of the Cartan group (i.e., characters stable under the Weyl group) \cite[Theorem 6.1]{Altug1}. This is indicative of how much more difficult the Poisson summation formula is in this context; rather than just the value of the Fourier transform at $\xi=0$ (which, as discussed, corresponds to the trace of the trivial representation), the point $\xi=0$ has a much richer contribution, including all the terms that the nonstandard nature of the space $\mathcal S_{r}((N,\psi)\backslash G/(N,\psi))$ is responsible for. On the spectral side of the nonstandard Kuznetsov formula, for the Eisenstein series unitarily induced from quadratic characters, it is reasonable to expect (and this will be made rigorous in my work in preparation with C.\ Wan, developing a version of the diagram \eqref{STFKTF} for $\GL_n$) that the double pole of the adjoint $L$-function at those points -- as opposed to its simple pole at other points of the continuous spectrum -- produces a discrete contribution to the nonstandard Kuznetsov formula.

Altu\u{g} goes on to insert the standard $L$-function into the trace formula in \cite{Altug3} (note that here we cannot keep working with $\SL_2$; we need to work with $G=\PGL_2$ or $\GL_2$, as Altu\u{g} does), restricted to holomorphic modular forms (by appropriate choice of test function at $\infty$). Note that the standard $L$-function of a cuspidal representation is entire, which at the level of $L$-groups reflects the fact that the stabilizer of a nonzero point in the standard representation belongs to a proper parabolic subgroup. Without taking the entire continuation as given, Beyond Endoscopy would like to study the residue 
\[ \sum_\pi \Res_{s=1} L(s,\pi),\]
with the sum ranging over modular forms of fixed weight $k$, and show that it is zero. This is what Altu\u{g} achieves by trace formula methods (for $k\ge 3$) showing, at the same time, how these methods can be used to obtain the holomorphic continuation of such an $L$-function beyond the domain of convergence.

\subsection{Comparison between the Jacquet and Kuznetsov trace formulas} \label{ss:Jacquet}

A global comparison analogous to \eqref{STFKTF} was established in \cite{SaBE2}, where the trace formula is replaced by the relative trace formula of Jacquet for Waldspurger periods, corresponding to a torus $T=\Res_{E/F}\Gm /\Gm \subset \PGL_2$, and $r=\Std \oplus \Std \otimes\eta$, where $\eta$ is the Galois character associated to the quadratic extension $E/F$.

\begin{equation}\label{RTFKTF}
\xymatrix{ \mathcal S(T\backslash G/T) \ar[rr]^{\mathcal T^{-1}}\ar[dr]_{\RTF} & & \ar[dl]^{\KTF}\mathcal S_{L(\Std \oplus \Std\otimes\eta, \frac{1}{2})}((N,\psi)\backslash G/(N,\psi))\\
& \mathbb C. &}
\end{equation}

In this case, there is no problem defining the functional $\RTF$ (as in \cite{Jacquet-Waldspurger}), while the functional $\KTF$ is beyond the domain of convergence, since it contains special values of the $L$-function associated to $r=\Std \oplus \Std\otimes\eta$ at the center of symmetry. The transfer operator $\mathcal T$, here, is given by the composition of $2$ Fourier transforms; more precisely, orbital integrals on both sides are parametrized by the affine line, and  $\mathcal T$ is given by a ``Kloosterman convolution:''
\begin{multline*} \mathcal Tf (u) =\eta(u) \cdot (\psi(\bullet)d\bullet)\star (\eta(\bullet)\psi(\bullet)d\bullet) \star f(u)\\ = \eta(u) \iint f\left(\frac{u}{xy}\right) \psi(x+y) \eta(x)  dx dy.\end{multline*}

Therefore, the statement that \eqref{RTFKTF} commutes can also be seen as a Poisson summation formula. Let us collect a few takeaways from the proof of this statement in \cite{SaBE2}.

\begin{enumerate}
 \item While the expression of $\KTF$ with nonstandard test measures corresponding to $L(\Std \oplus \Std\otimes\eta, \frac{1}{2})$ diverges, we can deform these test measures to correspond to $L(\Std \oplus \Std\otimes\eta, \frac{1}{2}+s)$, obtaining convergent expressions for $\Re(s)\gg 0$. There is nothing surprising here -- this just corresponds, classically, to feeding a series of Poincar\'e series into the Kuznetsov formula, corresponding to the Dirichlet series of the $L$-function -- but the same ``deformation'' can be applied to the transfer operator and the orbital integrals on the left hand side, obtaining a family of spaces $ \mathcal S_s(T\backslash G/T) $ of test measures which have nothing to do, a priori, with relative trace formulas. The natural extension of the functional $\RTF$ is evidently meromorphic on the entire plane, and the Poisson summation formula can be proven for $\Re(s)\gg 0$, leading to the meromorphic continuation of the right hand side. 
 
 \item One can encode the functional equation of the pertinent $L$-function into a ``Hankel transform'' on the space of test measures, as we saw in \S~\ref{sssHankel}:
 \[ L(\Std \oplus \Std\otimes\eta, \frac{1}{2}+s) \xrightarrow\sim L(\Std \oplus \Std\otimes\eta, \frac{1}{2}-s).\]
 This allows for a control of the functional $\KTF$ in the non-convergent region, and, by an application of the Phragm\"en--Lindel\"of principle, to the full spectral decomposition of the Kuznetsov formula in this region, and a new proof of the functional equation of the $L$-function. 
\end{enumerate}

The fundamental question here is if such ``transfer operators'' can be used to deal, similarly, with other $L$-functions inserted into the Kuznetsov formula.

\subsection{Functional equation of the standard $L$-function} \label{ss:Herman}

A Beyond-Endo\-scopy proof of the functional equation of an $L$-function had already appeared earlier in a paper of P.\ E.\ Herman \cite{Herman}, which has been little understood. The $L$-function here is the standard $L$-function for $\GL_2$, and Herman restricts himself to holomorphic cusp forms, by using the classical Petersson formula, instead of the full Kuznetsov formula. His proof is a masterful set of manipulations, which appear as rabbits pulled out of a hat, and lead to the following close analog of the Voronoi summation formula:
\begin{multline*} \sum_{f} \overline{a_l(f)} \sum_{n\geq 1} a_n(f)g(n) \\ = \sum_{f} \overline{a_{l}(f)}\left[ \frac{2\pi i^k\eta(f)}{\sqrt{D} }\sum_{n} a_n(f_D) \int_0^\infty g(x)J_{k-1}(\frac{4\pi\sqrt{nx}}{\sqrt{D}})dx\right].\end{multline*}
Here, $f$ ranges over an orthonormal basis of holomorphic modular forms of weight $k$, with some level $D$ and nebentypus $\chi$, and $g\in C_c^\infty(\mathbb R^+)$ is a smoothing function. (The integer $D$ is square-free, $\chi$ is a primitive Dirichlet character modulo $D$, and $k\ge 2$ is even; I point to the original article for the rest of the notation.)

It would be interesting to explore the connections between Herman's arguments and the Hankel transform $\mathcal H_{\Std}$  of Jacquet that we saw in \eqref{Hankel-Jacquet}. 
My expectation is that the analytic manipulations  in the proof of the formula above \cite[Theorem 3.1]{Herman} prove a Poisson summation formula for this Hankel transform, i.e., the statement that the following diagram commutes
\begin{equation}\label{HankelKTF}
\xymatrix{ \mathcal S_{L(\Std , \frac{1}{2}-s)}((N,\psi)\backslash G/(N,\psi)) \ar[rr]^{\mathcal H_{\Std}}\ar[dr]_{\KTF} & & \ar[dl]^{\KTF}\mathcal S_{L(\Std , \frac{1}{2}+s)}((N,\psi)\backslash G/(N,\psi))\\
& \mathbb C,&}
\end{equation}
of course under some simplifying assumptions (holomorphic modular forms, field of rationals, etc.). His manipulations definitely have semblance to Jacquet's formula: for example, the Poisson summation formula is applied twice, with a clever modification of an exponential factor in between; this is similar to the two Fourier transforms, with the intermediate correction factor, of \eqref{Hankel-Jacquet}.

\subsection{Lifting from tori to $\GL_2$} \label{ss:Venkatesh}

We now return to the setting of \S~\ref{ss:Kuztori}, where, for $G=\SL_2$ and $T$ a 1-dimensional torus, for the natural map $\iota:{^LT}\to {^L\SL_2}$ of $L$-groups, we described a transfer operator 
\[\mathcal S_{L(\Ad,\eta |\bullet|)} (N,\psi\backslash \SL_2/N,\psi) \xrightarrow{\mathcal T_\iota}\mathcal S(T),\] 
given by the formula \eqref{transfer-Venkatesh}.

Venkatesh's thesis can be understood as a statement that entails the commutativity of a diagram 
\begin{equation}\label{toriKTF}
\xymatrix{ \mathcal S_{L(\Ad,\eta |\bullet|)}((N,\psi)\backslash G/(N,\psi)) \ar[rr]^{\mathcal T_\iota}\ar[dr]_{\Res \KTF} & & \ar[dl]^{\TF}\mathcal S(T)\\
& \mathbb C,&}
\end{equation}
where ``$\Res \KTF$'' is supposed to contain only those $L$-packets of automorphic representations $\pi$ of $\SL_2$ for which $L(s,\Ad\otimes\eta, \pi)$ has a pole at $s=1$. 

To make sense of this, we enlarge our group to $\tilde G= \Gm\times \SL_2$. There is a certain nonstandard space of test measures $\mathcal S_{L(\Sym^2, 1)} (N,\psi\backslash \tilde G/N,\psi)$, constructed in \cite[\S~8.3]{SaTransfer2}, together with a $(\Gm,\eta)$-equivariant surjection
\[\mathcal S_{L(\Sym^2, 1)} (N,\psi\backslash \tilde G/N,\psi)\xrightarrow{p} \mathcal S_{L(\Ad,\eta |\bullet|)}((N,\psi)\backslash G/(N,\psi)).\]
As the notation suggests, the space $\mathcal S_{L(\Sym^2, 1)} (N,\psi\backslash \tilde G/N,\psi)$ contains the image $f_{L(\Sym^2, 1)}$ of the generating series of Hecke elements for $L(1,\Sym^2, \pi)$, where $\Sym^2$ is the symmetric square representation, which factors 
\[ \GL_2 \to {^L G} = \Gm \times \PGL_2 \to \GL_3.\]
If we use coordinates $(z,\zeta)$ for $N\backslash \tilde G \sslash N \simeq \Gm \times N\backslash \SL_2 \sslash N$, 
we can identify elements of $\mathcal S_{L(\Sym^2, 1)} (N,\psi\backslash \tilde G/N,\psi)$ with measures in the variables $(z,\zeta)$, so that the basic vector of $\mathcal S_{L(\Ad,\eta |\bullet|)}((N,\psi)\backslash G/(N,\psi))$ is the pushforward of the product $\eta(z) f_{L(\Sym^2, 1)}(z,\zeta)$ via the map $(z,\zeta)\mapsto \zeta$. More generally, we will consider the pushforward of $\eta(z)|z|^s f_{L(\Sym^2, 1)}(c,z)$,
which we will denote by $f_{\eta|\bullet|^s}$. We would then like to calculate 
\[ \Res_{s=1} \KTF(f_{\eta|\bullet|^s}).\]

We now write the composition of the transfer operator $\mathcal T_\iota$ with the projection $p$ as follows: write $f(z,\zeta) = \Phi(z,\zeta) \, d^\times z \, d\zeta$, then  
\begin{multline}\label{transform}\frac{\mathcal T_\iota \circ p (f)(t)}{dt} =  \lambda(\eta,\psi) \int_r \int_y \Phi\left( r, \frac{x(t)}{y} \right) |y|^{-1}  \eta(yrx(t)) \psi(y) dy d^\times r  \\
=  \lambda(\eta,\psi) \lim_{s\to 1} \int_r \eta(r) |r|^s\int_y \Phi\left( \frac{y r}{x(t)} , \frac{x(t)}{y} \right)|y|^{-1}   \psi(y) dy  \, d^\times r.\end{multline}

The inner integral of \eqref{transform} is a Fourier transform, and suggests applying a Poisson summation formula, globally. Denoting 
\[\hat f(r,x) = \int_y \Phi\left( \frac{y r}{x} , \frac{x}{y} \right)|y|^{-1}   \psi(y) dy  =  \int_u \Phi\left( ur , u^{-1} \right) |u|^{-1}   \psi(x u) du,\] the first step that Venkatesh performs is a Poisson summation formula which can be roughly thought of as an equality of the form 
\begin{align}  \nonumber \KTF(f_{\eta|\bullet|^s})  \approx \sum_{\zeta\in \QQ} \frac{f_{\eta|\bullet|^s}}{d\zeta}(\zeta) &=  \int_{\QQ^\times\backslash\mathbb A^\times} \eta(r)|r|^s \sum_{(n,c)\in (\QQ^\times \times \QQ)}  \Phi(rn,\frac{c}{n})  d^\times r \\ \label{KTF-Poisson}
&\approx \int_{\QQ^\times\backslash\mathbb A^\times} \eta(r)|r|^s \sum_{(\nu,c)\in \QQ^2}  \hat f(rc,\nu)  d^\times r,
\end{align}
where we have used the fact that the Fourier transform of $n\mapsto \Phi(rn, \frac{c}{n})|n|^{-1}$, with dual variable denoted $\frac{\nu}{c}$, is 
\[ \int_n \Phi(rn, \frac{c}{n}) |n|^{-1}\psi(n \frac{\nu}{c}) dn =  \int_y \Phi(\frac{rc y }{\nu}, \frac{\nu}{y}) |y|^{-1} \psi(y) dy =  \hat f(rc, \nu).\]
I warn the reader that in the thesis of Venkatesh, the local transfer operators were not available, and this is not the way that his theorem is presented; I point to \cite[\S~10.3]{SaTransfer2} for a slightly more expanded discussion of the relation.

The outer integral in \eqref{transform} is already embedded in \eqref{KTF-Poisson}. Classically, it corresponds to a Dirichlet series in the parameter $c$, with the parameter $\nu$ fixed. Studying this Dirichlet series is the second main step  of Venkatesh, where it is denoted by $Z(s-1)$ \cite[\S 4.5.2]{Venkatesh}. 
As Venkatesh observes, this Dirichlet series exhibits a lot of cancellation and does not have a pole at $s=0$, unless $\nu$ in the image of the map $T(\mathbb Q)\to T\sslash W (\mathbb Q)$. Locally, however, there is no pole, and the evaluation at $s=0$ corresponds to the integral of $\hat f$ against $\eta$ in the variable $c$. The reasons for these statements are essentially local: Fixing the parameter $\nu$, it turns out that the function $c\mapsto \eta_\nu(c) |c|^{-1} \hat f(c,\nu)$ is a smooth function around $c=0$; here, $\eta_\nu$ is the quadratic (or trivial) character corresponding to the splitting field of the characteristic polynomial $x^2-\nu x +1$. Moreover, for the purposes of analyzing the poles of the global zeta integral, this function at almost every non-Archimedean place $v$ can be approximated by the characteristic function of $1_{\mathfrak o_v}$, see \cite[Proposition 5]{Venkatesh}.  Thus, fixing $\nu\in \mathbb Q$, the last integral of \eqref{KTF-Poisson} looks like a Tate zeta integral, evaluated at the character $|\bullet|^{1+s}\eta\cdot \eta_\nu$, which explains why it only has a pole when $\eta_\nu = \eta$. 

Thus, using \eqref{transform}, and the fact that the Euler product of factors $\lambda(\eta,\psi)$ is $1$, we can write 
\[ \Res_{s=1} \KTF(f_{\eta|\bullet|^s}) = \sum_{\gamma\in T(\QQ)} \frac{\mathcal T_\iota \circ p (f)}{dt}(\gamma),\]
realizing the commutativity of \eqref{toriKTF}; of course, the ``trace formula'' for the torus is simply the sum of the test function over its rational points.

Finally, using the action of the Hecke algebra and the fundamental lemma  \eqref{FL-torus} for the transfer operator, we can isolate representations in the equality above, and prove functoriality for the morphism $\iota$ of $L$-groups \cite[Theorem 1]{Venkatesh}.

\subsection{General discussion}

What does the experience collected so far add to Langlands' original vision? The results to date are still too limited to pass verdict on the prospects of the program, but some structure has started emerging, which allows for cautious optimism. Based on this experience, I would like to end this article by ruminating on how we can move further.

\begin{enumerate}
 \item To begin with, since it is more manageable (and, geometrically, more natural) to insert $L$-functions, rather than their logarithmic derivatives, into trace formulas, one needs to rethink what the residues of this type of ``$r$-trace formula'' (where $r$ denotes the representation of the $L$-group whose associated $L$-function is to be inserted) should be compared with. If poles of $L$-functions detect functorial lifts from smaller groups, the residues (at least, when the poles are simple) are also themselves special values of $L$-functions, which should be inserted into the trace formula for the smaller group. (In that sense, the thesis of Venkatesh discussed in \S~\ref{ss:Venkatesh} is special, because the adjoint $L$-values appearing as denominators on the spectral side of the Kuznetsov formula essentially cancel out the residues.) A thoughtful strategy of how to do this systematically needs to be formulated.
 
 \item Related to the previous item, but already present in Langlands' original approach, is the question of ``dimension data'' for subgroups of the $L$-group. The question is whether the multiplicities of poles of $L$-functions are sufficient to detect the ``Zariski closure of a Langlands parameter.'' Although the notion of global Langlands parameters is conjectural, the question makes sense, thinking of poles as manifestations of the trivial representation, and asking questions about its multiplicity, when representations $r$ of the $L$-group are restricted to various closed subgroups. As was demonstrated by J.\ An, J.-K.\ Yu, and J.\ Yu \cite{AYY}, dimension data determine the conjugacy class of a subgroup up to finitely many possibilities. Langlands coined the term ``hadronic'' for those automorphic representations (of Ramanujan type) whose parameter does not land in a proper subgroup (and, hence, whose $L$-functions, for $r$ irreducible nontrivial, should have no poles). At the end of the day, the idea goes, there should be a way to break down the stable trace formula of a group into ``hadronic pieces'' for $G$ and smaller (dual) groups. 
 
 \item The point of view on Altu\u{g}'s work presented in \S~\ref{ss:Altug} suggests that, when removing non-Ramanujan and other ``special'' contributions from the stable trace formula, one ends up with the Kuznetsov formula. It is therefore unclear if anything is to be gained by working with the Arthur--Selberg trace formula, rather than directly with the Kuznetsov formula, as suggested early in the day by Sarnak \cite{Sarnak}. This is quite natural from the point of view of the geometric Langlands program, as well, where the Whittaker model has a very special place, essentially corresponding to the structure sheaf on the spectral side of the correspondence.
 
 \item 
As we have seen, to date there has only been one successful global ``Beyond Endoscopy'' comparison of (relative) trace formulas with \emph{unequal} underlying dual groups -- namely, the one in Venkatesh's thesis. Once the commutativity of a diagram such as \eqref{STFKTF} is established (expanding, slightly, upon Altu\u{g}'s work), the combination of the two will constitute a proof of functoriality from tori to $\GL_2$ in the way that Langlands originally envisioned -- using the Arthur--Selberg trace formula.

\item 
As we have seen in this section, there always seems to be a local comparison between the global comparisons established, even if that was not realized or used by the original authors. The local ``transfer operators'' could eventually guide our quest for a global comparison. In the most optimistic (probably too optimistic) scenario, they will be given by a combination of abelian Fourier transforms and other, rather innocuous, factors, which ``in principle'' admit a Poisson summation formula. Understanding the local transfer operators in arbitrary rank is the single most important task, in order to create good foundations for the Beyond Endoscopy program, and to formulate conjectures analogous to the Fundamental Lemma, that will guide the further development of the program.

\item 
Our discussion of global methods involved not only transfer operators, but also the Hankel transforms which encode the functional equations of $L$-functions, that we encountered in \S~\ref{ss:Hankel}. Their appearance is very natural, once $L$-functions are inserted into the trace formula, since one wants to move beyond the domain of their convergence. Indeed, Hankel transforms play a role in \S~\ref{ss:Jacquet} in establishing a \emph{full} meromorphic continuation of the nonstandard (Kuznetsov) trace formula, avoiding difficult and ad hoc analytic arguments, and they are probably behind the ``Beyond Endoscopy'' proof of the functional equation discussed in \S~\ref{ss:Herman}. This supports the idea that the topics that we encountered in this article -- Langlands' strategy for proving functoriality, and the study of $L$-functions by means of Schwartz spaces and Fourier/Hankel transforms -- are very closely related, and will all have a role to play beyond endoscopy.

\end{enumerate}

\bibliographystyle{amsalpha}
\bibliography{biblio}

\end{document}